\documentclass[12pt]{amsart}
\usepackage{amsmath}
\usepackage{amsthm}
\usepackage{amssymb}
\usepackage{tikz-cd}
\usepackage{hyperref}

\makeatletter
\def\l@subsection{\@tocline{2}{0pt}{2.5pc}{5pc}{}}
\makeatother

\usepackage{url}
\usepackage{cite}
\usepackage{dsfont}
\usepackage{array}

\usepackage[all]{xy}
\usepackage{mathrsfs}
\usepackage{amscd}
\usepackage{latexsym}
\usepackage{epsfig}
\usepackage{graphicx}
\usepackage{amsfonts}
\usepackage{psfrag}
\usepackage{caption}
\usepackage{rotating}
\usepackage{mathtools}
\usepackage{etoolbox}

\usepackage{bm}

\newtheorem{theorem}{Theorem}[section]
\newtheorem{proposition}[theorem]{Proposition}
\newtheorem{corollary}[theorem]{Corollary}

\theoremstyle{definition}
\newtheorem{definition}[theorem]{Definition}

\theoremstyle{remark}
\newtheorem*{remark}{Remark}

\usepackage{chngcntr}
\counterwithin{equation}{section}

\DeclareMathOperator{\Ad}{Ad}
\DeclareMathOperator{\ad}{ad}
\DeclareMathOperator{\KK}{\mathbf{K}}
\DeclareMathOperator{\kk}{\mathbf{k}}

\DeclareMathOperator{\tr}{tr}

\DeclareMathOperator{\SU}{\mathbf{SU}}
\DeclareMathOperator{\su}{\mathfrak{su}}

\DeclareMathOperator{\Hom}{Hom}
\DeclareMathOperator{\End}{End}
\DeclareMathOperator{\Aut}{Aut}
\DeclareMathOperator{\Stab}{Stab}
\DeclareMathOperator{\NN}{N}

\DeclareMathOperator{\ZZ}{Z}
\DeclareMathOperator{\ee}{\mathrm{e}}

\DeclareMathOperator{\pr}{pr}
\DeclareMathOperator{\vol}{vol}

\DeclareMathOperator{\res}{res}

\DeclareMathOperator{\sgn}{sgn}
\DeclareMathOperator{\Id}{Id}
\DeclareMathOperator{\dimr}{\dim_{\mathbb{R}}}

\DeclareMathOperator{\diag}{diag}

\newcommand{\I}{\mathbf{i}}
\newcommand{\orb}{\mathcal{O}}
\newcommand{\lieg}{\mathfrak{g}}
\newcommand{\liegdual}{\mathfrak{g}^{*}}
\newcommand{\liet}{\mathfrak{t}}
\newcommand{\lietdual}{\mathfrak{t}^{*}}

\newcommand{\cpwc}{\mathfrak{t}_{\geq 0}}
\newcommand{\opwc}{\mathfrak{t}_{> 0}}
\newcommand{\diw}{\Lambda_{\geq 0}}
\newcommand{\proots}{\mathrm{R}_{+}}

\newcommand{\red}{\mathrm{red}}

\newcommand{\CC}{\mathbb{C}}
\newcommand{\cc}{\mathbf{c}}
\newcommand{\RR}{\mathbb{R}}
\newcommand{\FFF}{\mathbb{F}}
\newcommand{\ZZZ}{\mathbb{Z}}
\newcommand{\MM}{\mathcal{M}}
\newcommand{\veceta}{\vec{\eta}}
\newcommand{\vecxi}{\vec{\xi}}
\newcommand{\vecw}{\vec{w}}
\newcommand{\ws}{\mathfrak{s}}

\newcommand{\SSS}{\mathcal{S}}
\newcommand{\RRR}{\mathcal{R}}

\newcommand{\calK}{\mathcal{K}}

\title{Volume Formula  for $ N $-fold Reduced Products}

\author{Lisa Jeffrey}\thanks{The author is partially supported by an NSERC Discovery Grant.}
\address{Mathematics Department, University of Toronto, Toronto, ON M5S 2E4, Canada} 
\author{Jia Ji}\thanks{The author is partially supported by the University of Toronto.}
\address{Mathematics Department, University of Toronto, Toronto, ON  M5S 2E4, Canada}

\begin{document}


\maketitle


\begin{abstract}
Let $ G $ be a semisimple compact connected Lie group.
An $ N $-fold reduced product of $ G $ is the symplectic
quotient of the Hamiltonian system of the Cartesian product
of $ N $ coadjoint orbits of $ G $ under diagonal coadjoint action of $ G $. Under appropriate
assumptions, it is a symplectic orbifold. Using the technique
of nonabelian localization and the residue formula of Jeffrey
and Kirwan, we investigate the symplectic volume 
 of an $ N $-fold
reduced product of $ G $. Suzuki and Takakura 
gave a volume formula for the  $ N $-fold reduced product of $ \mathbf{SU}(3) $ in \cite{ST08}
by using geometric quantization and
the Riemann-Roch formula.  We compare our volume formula with theirs
and prove that our volume formula
agrees with theirs in the case of triple reduced products
of $ \mathbf{SU}(3) $.
\end{abstract}






\tableofcontents










\section{Introduction}


In this article, we study objects called $ N $-fold reduced products. The formal definitions
will be given in later chapters. Roughly speaking, these objects are reduced spaces obtained
from certain Hamiltonian systems.

Given a semisimple compact connected Lie group $ G $, we may consider
its adjoint orbits in its Lie algebra $ \lieg $.
Notice that in this article we will fix a $ G $-invariant inner product on $ \lieg $ and hence we can
identify adjoint orbits with coadjoint orbits, which are naturally symplectic manifolds.
If we are given $ N $ such adjoint orbits,
we can form their Cartesian product denoted by $ M $. The group acts on this product diagonally through adjoint
action. This is a Hamiltonian system with moment map $ \mu_{G}: M \to \lieg $.
An $ N $-fold
reduced product is the reduced space $ \mu_{G}^{-1}(0) / G $ of this Hamiltonian system.

In general, the geometry of these reduced products are very complicated and thus difficult to study.
For example, they are in general not smooth manifolds. Even with the assumptions that we are going to
make, they are still in general only orbifolds. Roughly speaking, an orbifold can be thought as
a space which is almost smooth with mild singularities (see \cite{Satake} and \cite{Duistermaat}).
By the Marsden-Weinstein reduction theorem, such a reduced product naturally carries a symplectic
structure and we are interested in understanding its symplectic geometry. 

One of the most important
symplectic invariants is the symplectic volume and this is the topic we will focus on the most in this
article. We will use the technique of nonabelian localization and the residue formula, developed by
Jeffrey and Kirwan (\cite{JK95}, \cite{JK3}), to study the symplectic volume and also the intersection
pairings of these reduced products. 

In 2008, Suzuki and Takakura studied the symplectic volume of
$ N $-fold reduced products of $ G = \SU(3) $ in their paper \cite{ST08} via
Riemann-Roch. Our volume formula generalizes their volume formula in the sense that their volume
formula requires more restrictive inputs. Furthermore, in the case of $ N = 3 $,
i.e., the case of triple reduced products of $ \SU(3) $, we have proved that
up to normalization constants, our volume formula  completely agrees with theirs.

In a later paper \cite{ST09}, Suzuki and Takakura generalize their results 
on symplectic volumes of $N$-fold 
reduced products of $G =\SU(3)$ to the case where $G$ is a general 
compact Lie group.  We are able to obtain results for the case 
of general $G$ as well. Our results appear in \S 
\ref{ch:generalizations} below. Our methods may be  
adapted to generalizing from volumes to intersection pairings.
We address that subject in \cite{JJ19}.

$ N $-fold reduced products appear in the literature in other guises as well.
For example,  $ N $-fold reduced products of $ G = \SU(2) $ can be identified with
moduli spaces of polygons in $ \RR^{3} $ with prescribed lengths of edges, which have been
studied by, for example, Hausmann and Knutson \cite{HK}, Kamiyama and Tezuka \cite{KT99}.
By Jeffrey \cite{Jeffrey} (Theorem 6.6 in that paper), $ N $-fold reduced products of $ G $ can also be identified with moduli spaces of 
flat
$ G $-connections on a genus $ 0 $ surface with $ N $ boundary components and the corresponding holonomies conjugate
to $ \exp(\xi_{i}) $ for prescribed $ \xi_{i} \in \lieg $ (provided that the $ \xi_{i} $'s are sufficiently
small).

The organization of this article is as follows.

In Section \ref{ch:trp}, we investigate triple reduced products of $ \SU(3) $, namely, the case of $ G = \SU(3) $ and $ N = 3 $.
Along the way, we also introduce the notations that can be easily generalized in the later chapters.
After briefly reviewing the machinery of equivariant cohomology, we describe the method of nonabelian
localization. 
In Section 3, we review  and the residue formula of Jeffrey and Kirwan (\cite{JK95}, \cite{JK3}). Then we
apply these techniques to derive our volume formula. In 
Section 4, 
we compare our volume formula with the volume formula of Suzuki and Takakura (\cite{ST08}) and conclude
this chapter with a proof that our volume formula matches
with theirs in the case of triple reduced products of $ \SU(3) $.

Finally in Section \ref{ch:generalizations} we first generalize our volume formula of triple reduced products
of $ \SU(3) $ to the case of $ N $-fold reduced products of $ \SU(3) $. Then we further generalize it
to the case of $ N $-fold reduced products of a general semisimple compact connected Lie group $ G $.

In a later article \cite{JJ19}, by applying the residue formula of Jeffrey and Kirwan (\cite{JK95},\cite{JK3}), we compute the intersection pairings of $ N $-fold
reduced products.


\section{Background} \label{ch:trp}

In this section, we set the stage by first reviewing some basic facts about $ \SU(3) $
and then constructing triple reduced products from adjoint orbits of $ \SU(3) $.
Along the way we set the notations that will be used throughout this article.

\subsection{Basic Facts about $ \SU(3) $} 
We review the following basic material.
\begin{itemize}
\item Let $ G = \SU(3) $,   the collection of all invertible
complex $ 3 \times 3 $ matrices $ A $ such that $ A^{*} = A^{-1} $
and $ \det(A) = 1 $, where $ A^{*} $ denotes the conjugate transpose of $ A $.
\item Let $ \lieg = \su(3) $, the Lie algebra of $ \SU(3) $.
The vector space  $ \lieg $ is the collection of all complex $ 3 \times 3 $ matrices $ X $
such that $ X^{*} = -X $ and $ \tr(X) = 0 $. Notice that $ \lieg $ is
a vector space over $ \RR $ and $ \dimr(\lieg) = 8 $.
\item Let $ T $ be the standard maximal torus of $ G $. In other words, $ T $ is the
collection of all $ 3 \times 3 $ diagonal matrices $ \diag(e^{\I \theta_{1}}, 
e^{\I \theta_{2}}, e^{-\I (\theta_{1} + \theta_{2})}) $ such that $ \theta_{1}, 
\theta_{2} $ are real.
\item Let $ \liet $ be the
Lie algebra of $ T $. Thus $ \liet $ is the collection of all 
$ 3 \times 3 $ diagonal
matrices $ \diag(\I \theta_{1}, \I \theta_{2}, -\I (\theta_{1} + \theta_{2})) $ such that
$ \theta_{1}, \theta_{2} $ are real.
\item Let $ \liegdual := \Hom_{\RR}(\lieg, \RR) $ be the dual vector space of $ \lieg $, and
let $ \lietdual := \Hom_{\RR}(\liet, \RR) $ be the dual vector space of $ \liet $.
\item If $ V $ is a vector space over a field $ \FFF $,
let $ \left\langle \xi, X \right\rangle := \xi(X) \in \FFF $ denote the natural 
pairing between a covector $ \xi \in V^{*} := \Hom_{\FFF}(V, \FFF) $ and a
vector $ X \in V $. In this article, $ \FFF $ is either $ \RR $ or $ \CC $,
depending on the context.
\item For all $ g \in G $, let $ \cc(g): G \to G $ denote the map $ x \mapsto g x g^{-1} $.
\item For all $ g \in G $, let $ \Ad(g): \lieg \to \lieg $ denote the differential of
$ \cc(g) $ at the identity $ e \in G $. Thus, $ \Ad: G \to \Aut(\lieg) $ is the
adjoint representation of $ G $ on $ \lieg $. The differential of $ \Ad $ at
the identity $ e \in G $ is denoted by $ \ad $. Thus,
$ \ad: \lieg \to \End(\lieg) $ is the adjoint representation of $ \lieg $ on $ \lieg $.

\item For all $ g \in G $, let $ \KK(g): \liegdual \to \liegdual $ be defined by:
\begin{equation}
\left\langle \KK(g) \xi, X \right\rangle = \left\langle \xi, \Ad(g^{-1}) X \right\rangle 
\end{equation}
for all $ \xi \in \liegdual $, $ X \in \lieg $. Thus, $ \KK: G \to \Aut(\liegdual) $
is the coadjoint representation of $ G $ on $ \liegdual $. The differential of $ \KK $
at the identity $ e \in G $ is denoted by $ \kk $. Thus,
$ \kk: \lieg \to \End(\liegdual) $ is the coadjoint representation of $ \lieg $
on $ \liegdual $. Notice that for all $ X \in \lieg $,
\begin{equation}
\kk(X) = -\ad(X)^{*}.
\end{equation}
\item For all $ X, Y \in \lieg $, we define
\begin{equation}
(X, Y) := -\tr(XY).
\end{equation}
Then $ (\cdot, \cdot) $ is an $ \Ad(G) $-invariant (or briefly, $ G $-invariant) 
inner product on $ \lieg $.
We will use this inner product to identify $ \liegdual $ with $ \lieg $ through
the standard identification between $ X \in \lieg $ and $ (X, \cdot) \in \liegdual $
for all $ X \in \lieg $. Since $ \liet $ is naturally a subspace of $ \lieg $
through the inclusion $ T \subset G $, the above identification induces an identification
between $ \liet $ and $ \lietdual $. It is in this sense that we write $ \lietdual \subset \liegdual $.

\item Suppose $ \xi \in \liegdual $ and $ X \in \lieg $ satisfy
\begin{equation}
\left\langle \xi, Y \right\rangle = (X, Y)
\end{equation}
for all $ Y \in \lieg $. Namely, $ \xi $ corresponds to $ X $ under the identification
through the inner product. For all $ g \in G $, $ Y \in \lieg $,
\begin{align}
\left\langle \KK(g) \xi, Y \right\rangle &= \left\langle \xi, \Ad(g^{-1}) Y \right\rangle \\
&= (X, \Ad(g^{-1}) Y) \\
&= (\Ad(g) X, Y).
\end{align}
Thus, $ \KK(g) \xi $ corresponds to $ \Ad(g) X $ for all $ g \in G $ under the identification through
the inner product. It is in this sense that we say the coadjoint action of $ G $ on $ \liegdual $
corresponds to the adjoint action of $ G $ on $ \lieg $, and it is in this sense that
we identify coadjoint orbits with adjoint orbits. In this article, we will try to focus on
the adjoint version of the theory since it is often easier to carry out computations with
elements in $ \lieg $, which are matrices, than with elements in $ \liegdual $.

\item The Weyl group $ W $ is defined by $ \NN(T) / T $, where $ \NN(T) $ denotes the normalizer of
$ T $ in $ G $. Suppose $ g \in \NN(T) $. Then $ g $ is a representative of one Weyl group element
$ w \in W $. This element $ w $ acts on $ T $ through
\begin{equation}
w \cdot t = g t g^{-1}
\end{equation}
for all $ t \in T $. The element $ w $ acts on $ \liet $ through
\begin{equation}
w \cdot X = \Ad(g) X
\end{equation}
for all $ X \in \liet $.

\item Since $ G = \SU(3) $, $ W $ is isomorphic to the permutation group $ \mathfrak{S}_{3} $
of $ 3 $ letters. More precisely, let
\begin{equation}
W = \left\lbrace \ws_{0}, \ws_{1}, \ws_{2}, \ws_{3}, \ws_{4}, \ws_{5} \right\rbrace 
\end{equation}
where
\begin{align}
\ws_{0} \  &\text{corresponds to} \  \Id \in \mathfrak{S}_{3}, \label{eq:wsbegin} \\
\ws_{1} \  &\text{corresponds to} \  (1 \  2) \in \mathfrak{S}_{3}, \\
\ws_{2} \  &\text{corresponds to} \  (1 \  2 \  3) \in \mathfrak{S}_{3}, \\
\ws_{3} \  &\text{corresponds to} \  (1 \  3) \in \mathfrak{S}_{3}, \\
\ws_{4} \  &\text{corresponds to} \  (1 \  3 \  2) \in \mathfrak{S}_{3}, \\
\ws_{5} \  &\text{corresponds to} \  (2 \  3) \in \mathfrak{S}_{3}. \label{eq:wsend}
\end{align}
The signature $ \sgn(w) $ of a Weyl group element $ w \in W $ is defined as
the signature of its corresponding element $ \sigma \in \mathfrak{S}_{3} $, i.e.
\begin{equation}
\sgn(w) := \sgn(\sigma).
\end{equation}
Notice that we have chosen the subscripts $ j $ in $ \ws_{j} $ so that
\begin{equation}
\sgn(\ws_{j}) = (-1)^{j}.
\end{equation}
This will be convenient for later use.
We shall introduce a chosen basis for $\liet$ below, and  
express $W$ in terms of this basis.  In   (\ref{weylmatrix}) we are able to 
write $W$ as $2 \times 2$ matrices in terms of this basis.  

\item The Weyl group elements act on diagonal matrices by permuting the diagonal
entries. More precisely, if $ \sigma \in \mathfrak{S}_{3} $, then
$ \sigma $ corresponds to a Weyl group element $ w \in W $ and
\begin{equation}
w \cdot \diag(a_{1}, a_{2}, a_{3}) = \diag(a_{\sigma^{-1}(1)}, a_{\sigma^{-1}(2)}, a_{\sigma^{-1}(3)}).
\end{equation}
For example,
\begin{equation}
\ws_{2} \cdot \diag(a_{1}, a_{2}, a_{3}) = \diag(a_{3}, a_{1}, a_{2}).
\end{equation}
If $ \sigma \in \mathfrak{S}_{3} $ corresponds to a Weyl group
element $ w \in W $, we define
\begin{equation}
\sigma \cdot X := w \cdot X
\end{equation}
for all $ X \in \liet $.

\item Since here $ G = \SU(3) $, we can identify the Weyl group $ W $ with $ \mathfrak{S}_{3} $
in the above way.

\item Since the inner product $ (\cdot, \cdot) $ is $ \Ad(G) $-invariant on $ \lieg $,
it is $ W $-invariant on $ \liet $.

\item Let $ \mathcal{L} $ denote the integral lattice in $ \liet $, that is,
\begin{equation}
\mathcal{L} := \left\lbrace H \in \liet \  : \  e^{2 \pi H} = I \right\rbrace .
\end{equation}
Thus,
\begin{equation}
\mathcal{L} = \left\lbrace m_{1} H_{1} + m_{2} H_{2} \  : \  m_{1}, m_{2} \in \ZZZ \right\rbrace 
\end{equation}
where
\begin{align}
H_{1} &= \diag(\I, -\I, 0), \\
H_{2} &= \diag(0, \I, -\I).
\end{align}

\item A weight is an element of $ \lietdual $ such that it takes integer values on $ \mathcal{L} $.
We will use the identification through the inner product $ (\cdot, \cdot) $ to regard weights as elements
in $ \liet $.

Let $ \lambda_{1}, \lambda_{2} $ be elements in $ \liet $ such that
\begin{equation}
(\lambda_{i}, H_{j}) = \delta_{ij}
\end{equation}
for all $ i, j \in \left\lbrace 1, 2 \right\rbrace  $.
Then, we obtain
\begin{align}
\lambda_{1} &= \diag(\frac{2}{3} \I, -\frac{1}{3} \I, -\frac{1}{3} \I), \\
\lambda_{2} &= \diag(\frac{1}{3} \I, \frac{1}{3} \I, -\frac{2}{3} \I).
\end{align}

\item A weight $ \lambda \in \liet $ is called dominant if and only if $ (\lambda, H_{j}) \geq 0 $
for all $ j $. A weight $ \lambda \in \liet $ is called integral if and only if
$ (\lambda, H_{j}) \in \ZZZ $ for all $ j $. Thus,
\begin{equation}
\diw := \left\lbrace n_{1} \lambda_{1} + n_{2} \lambda_{2} \  : \  n_{1}, n_{2}
\  \text{are both nonnegative integers} \right\rbrace  
\end{equation}
is the collection of dominant integral weights, and
\begin{align}
\cpwc &:= \left\lbrace c_{1} \lambda_{1} + c_{2} \lambda_{2} \  : \  c_{1}, c_{2} \geq 0 \right\rbrace, \\
\opwc &:= \left\lbrace c_{1} \lambda_{1} + c_{2} \lambda_{2} \  : \  c_{1}, c_{2} > 0 \right\rbrace 
\end{align}
are the closed positive Weyl chamber and the open positive Weyl chamber, respectively.

\item Let
\begin{align}
\alpha_{1} &:= 2 \lambda_{1} - \lambda_{2} = H_{1}, \\
\alpha_{2} &:= - \lambda_{1} + 2 \lambda_{2} = H_{2}
\end{align}
be the standard simple roots for $ G = \SU(3) $.
Following the convention of \cite{JK95}, let
\begin{equation}
\proots = \left\lbrace 2 \pi \alpha_{1}, 2 \pi \alpha_{2}, 2 \pi (\alpha_{1} + \alpha_{2}) \right\rbrace 
\end{equation}
denote the set of $ 2 \pi $-modified positive roots of $ \SU(3) $. From now on, if we talk about
roots, we mean the $ 2 \pi $-modified ones unless stated otherwise.
\end{itemize}
\subsection{Construction of Triple Reduced Products of $ \SU(3) $}
\subsubsection{Set-up}
Let $ \xi \in \lieg $. Let $ \orb_{\xi} $ denote the adjoint orbit through $ \xi $.
Under the identification through the inner product $ (\cdot, \cdot) $, it is equivalent
to consider either adjoint orbits or their corresponding coadjoint counterparts. In
this article we will mainly use the adjoint setting.

It is well known that every coadjoint orbit admits a natural symplectic form,
called the Kirillov-Kostant-Souriau form, or briefly the KKS form. In the adjoint
setting, it is defined as follows.

Suppose $ \orb $ is an adjoint orbit in $ \lieg $ and $ \xi \in \orb $. Then the tangent
space of $ \orb $ at $ \xi $, $ T_{\xi} \orb $, which is naturally a subspace of $ \lieg $,
is the following collection:
\begin{equation}
T_{\xi} \orb = \left\lbrace \ad(X) \xi \  : \  X \in \lieg \right\rbrace.
\end{equation}
Notice that $ \ad(X) \xi = [X, \xi] $ for all $ X \in \lieg $.
\begin{definition}
The KKS form $ \omega $ on $ \orb $ is defined by
\begin{equation}
\omega_{\xi}([X, \xi], [Y, \xi]) := (\xi, [X, Y])
\end{equation}
for all $ \xi \in \orb $, $ X, Y \in \lieg $.
\end{definition}

The KKS form $ \omega $ is a closed nondegenerate $ 2 $-form.
Equipped with $ \omega $, $ \orb $ becomes a compact symplectic manifold.
The group $ G $ acts on $ \orb $ by the adjoint action and this makes $ \orb $
a Hamiltonian $ G $-space with the inclusion map $ \mu_{\orb}: \orb \hookrightarrow \lieg $
as the moment map. In other words, for all $ X \in \lieg $,
\begin{equation}
d \mu_{\orb}^{X} = \iota_{X^{\sharp}} \omega
\end{equation}
where $ \mu_{\orb}^{X}(\xi) := (\mu_{\orb}(\xi), X) = (\xi, X) $ for all $ \xi \in \orb $ and $ X^{\sharp} $ is the fundamental vector field
on $ \orb $ generated by $ X $ and thus $ X^{\sharp}(\xi) = [X, \xi] \in T_{\xi} \orb $
for all $ \xi \in \orb $. Notice that
\begin{equation}
([X, Y], \xi) = (X, [Y, \xi])
\end{equation}
for all $ X, Y, \xi \in \lieg $.

Consider $ 3 $ points $ \xi_{1}, \xi_{2}, \xi_{3} $ in $ \lieg $. We then have $ 3 $
adjoint orbits $ \orb_{\xi_{1}}, \orb_{\xi_{2}}, \orb_{\xi_{3}} $. Let $ \omega_{i} $
denote the KKS form on $ \orb_{\xi_{i}} $.
Now consider $ M := \orb_{\xi_{1}} \times \orb_{\xi_{2}} \times \orb_{\xi_{3}} $.
Let $ \pr_{i}: M \to \orb_{\xi_{i}} $ be the standard projection onto the $ i $-th
factor. Then the form
\begin{equation}
\omega := \sum_{i = 1}^{3} \pr_{i}^{*} \omega_{i}
\end{equation}
is a symplectic form on $ M $.
\subsubsection{Triple reduced products}

Let $ G $ act on $ M $ by
\begin{equation}
g \cdot (\eta_{1}, \eta_{2}, \eta_{3}) := (\Ad(g) \eta_{1}, \Ad(g) \eta_{2}, \Ad(g) \eta_{3})
\end{equation}
for all $ g \in G $, $ \eta_{i} \in \orb_{\xi_{i}} $. Then this action
makes $ M $ a Hamiltonian $ G $-space with the moment map $ \mu_{G}: M \to \lieg $ such that
\begin{equation}
\mu_{G}(\eta_{1}, \eta_{2}, \eta_{3}) = \sum_{i = 1}^{3} \mu_{\orb_{\xi_{i}}}(\eta_{i})
= \sum_{i = 1}^{3} \eta_{i}
\end{equation}
for all $ \eta_{i} \in \orb_{\xi_{i}} $.

\subsubsection{Assumptions} \label{sss:assumptions}
In this article, we make the following assumptions about the points $ \xi_{1}, \xi_{2}, \xi_{3} $:
\begin{description}
\item[(A1)] $ M_{0} := \mu_{G}^{-1}(0) \neq \emptyset $ and $ 0 \in \lieg $ is a regular
value of $ \mu_{G} $.
\item[(A2)] $ \xi_{i} \in \cpwc \subset \lieg $ for all $ i $ and
each $ \orb_{\xi_{i}} $ is diffeomorphic to the homogeneous space $ G / T $.
\end{description}

\begin{remark}
The assumption (A1) ensures that   the stabilizer
$ \Stab_{G}(\veceta) $ is finite for each $ \veceta
= (\eta_{1}, \eta_{2}, \eta_{3}) \in M_{0} $. Every adjoint orbit $ \orb $
will intersect $ \cpwc \subset \lieg $ at exactly one point $ \xi $, so by only
considering $ \xi \in \cpwc $, we still obtain every possible orbit.
This explains the first part of (A2). The second part of (A2) says that
the orbits we will consider are nondegenerate, that is, they are of the highest dimension
possible. In fact, assuming (A2) is equivalent to assuming that
$ \xi_{i} \in \opwc \subset \lieg $ for all $ i $.
\end{remark}

Let $ M^{T} $ denote the set of fixed points in $ M $ under the action of $ T \subset G $.
We have the following.

\begin{proposition}\label{prop:TfixedPoints.N=3}
Let $ M = \orb_{\xi_{1}} \times \orb_{\xi_{2}} \times \orb_{\xi_{3}} $ be the Cartesian
product of $ N = 3 $ adjoint orbits of $ G = \SU(3) $, where the $ \xi_{i} $ satisfy
the assumptions (A1) and (A2). Then,
$ M^{T} $ is the discrete set
\begin{equation}
\left\lbrace (w_{1} \cdot \xi_{1}, w_{2} \cdot \xi_{2}, w_{3} \cdot \xi_{3}) \  : \  
w_{i} \in W \right\rbrace .
\end{equation}
Thus, $ \left| M^{T} \right| = \left| W \right| ^{3} $.
\end{proposition}

\begin{proof}
Note that
\begin{equation}
\orb_{\xi_{i}} \cap \liet = W \cdot \xi_{i}
\end{equation}
for all $ i $. The elements of $ \orb_{\xi_{i}} \cap \liet $ are precisely
those elements in $ \orb_{\xi_{i}} $ that are fixed by the adjoint action
of $ T \subset G $. Since $ \xi_{i} \in \opwc $, $ \orb_{\xi_{i}} \cap \liet $
is discrete and has the same cardinality as $ W $. The proposition follows
immediately.
\end{proof}

Since $ G = \SU(3) $, we actually have $ 6^{3}  $ isolated fixed points in
$ M $ under the action of $ T \subset G $.

\begin{definition}
The quotient space
\begin{equation}
M_{\red} := M_{0} / G
\end{equation}
is called a \emph{triple reduced product of $ G $}. In other words,
$ M_{\red} $ is the Marsden-Weinstein reduction of the Hamiltonian
$ G $-space $ (M, \omega, G, \mu_{G}) $.
Sometimes we may write $ M_{\red}(\vecxi) $ to emphasize the dependence on
the initial data $ \vecxi = (\xi_{1}, \xi_{2}, \xi_{3}) $.
\end{definition}

\begin{remark}
In general, $ M_{\red} $ is not a smooth
manifold. It belongs to a type of spaces called \emph{orbifolds} or
\emph{V-manifolds} (\cite{Satake}). Roughly speaking, an orbifold
is almost a smooth manifold except that it has some mild singularities.
At these singularities, it locally looks like $ U / \Gamma $ where
$ U $ is an open subset of $ \RR^{d} $ and $ \Gamma $ is a finite group
of linear automorphisms of $ U $. 

Fortunately, \cite{Satake} tells us
that on such spaces the de Rham theory and Poincar\'{e} duality work
basically the same way as on smooth manifolds (See also \cite{Duistermaat}).
However, to avoid such complication, we add the following assumption
about the $ \xi_{i} $:
\begin{description}
\item[(A3)] $ M_{\red}(\vecxi) $ is a smooth manifold.
\end{description}
This assumption is equivalent to the assumption that for each $ \veceta \in M_{0} $,
$ \Stab_{G}(\veceta) = \ZZ(G) $.
\end{remark}

\subsection{Equivariant Cohomology} \label{ss:eqc}

To study the symplectic volume $ \vol^{\SSS}(M_{\red}) $, we are basically looking at
the cohomological quantity
\begin{equation}
e^{\omega_{\red}} [M_{\red}] = \frac{1}{\I^{d / 2}} e^{\I \omega_{\red}} [M_{\red}],
\end{equation}
where $ [M_{\red}] $ denotes the fundamental class of $ M_{\red} $, which is picked up
by the orientation induced by the symplectic volume form on $ M_{\red} $.

This quantity can be computed using the technique called  nonabelian localization due
to Jeffrey and Kirwan (\cite{JK95}; see also \cite{JK3}) and in particular the residue formula
(Theorem 8.1 in \cite{JK95} and Theorem 3.1 in \cite{JK3}).
To state their results, we shall first review the machinery of equivariant cohomology.

In this article, we only consider cohomology groups over $ \CC $.
Let $ \calK $ be a compact connected Lie group. Let $ \MM $ be a $ \calK $-space.
The equivariant cohomology  of the $ \calK $-space $ \MM $ is
denoted $ H_{\calK}^{*}(\MM) $.
\begin{equation}
\MM \times_{\calK} E\calK := (\MM \times E\calK) / \calK,
\end{equation}

We will use the Cartan model (\cite{Cartan1,Cartan}; see also \cite{BGV}) to compute equivariant cohomology.
Let $ \mathfrak{k} $ denote the Lie algebra of $ \calK $.
A $ \calK $-equivariant differential form $ \alpha $ on $ \MM $ can be thought of as a $ \calK $-equivariant
polynomial map
\begin{equation}
\alpha: \mathfrak{k} \to \Omega^{*}(\MM).
\end{equation}
Let $ \Omega_{\calK}^{*}(\MM) $ denote the collection of all $ \calK $-equivariant differential forms
on $ \MM $. In other words,
\begin{equation}
\Omega_{\calK}^{*}(\MM) = (S(\mathfrak{k}^{*}) \otimes \Omega^{*}(\MM))^{\calK}.
\end{equation}

Now we return to our situation, that is, the situation of triple reduced products of $ G = \SU(3) $.

There is another important map called the pushforward map
\begin{equation}
\Pi^{G}_{*}: H^{*}_{G}(M) \to H^{*}_{G},
\end{equation}
which can be thought of as integration over $M$.
Here we have introduced the notation $ H^{*}_{G}$ to mean
$ H^{*}_{G}({\rm pt})$,
 and likewise $H^*_T$ for $H^*_T({\rm pt})$.
Similarly, when we are looking at the $ T $-action on $ M $, the corresponding pushforward map
is
\begin{equation}
\Pi^{T}_{*}: H^{*}_{T}(M) \to H^{*}_{T}.
\end{equation}
Usually, when the context is clear, we will denote both $ \Pi^{G}_{*} $ and $ \Pi^{T}_{*} $
by the integration symbol $ \int_{M} $ or simply $ \Pi_{*} $.

Recall that $ 0 $ is a regular value for the moment map $ \mu_{G} $.
By \cite{Kirwan}, the ring homomorphism
\begin{equation}
i_{0}^{*}: H^{*}_{G}(M) \to H^{*}_{G}(M_{0})
\end{equation}
is surjective. In addition, we have a canonical isomorphism
\begin{equation}
\pi_{0}^{*}: H^{*}(M_{\red}) \to H^{*}_{G}(M_{0})
\end{equation}
induced from the map
\begin{equation}
\pi_{0}: M_{0} \times_{G} EG \to M_{\red}.
\end{equation}
We have a a canonical map (the Kirwan map):
\begin{equation}
\kappa_{0} := (\pi_{0}^{*})^{-1} \circ i_{0}^{*}: H^{*}_{G}(M) \to H^{*}(M_{\red}).
\end{equation}

In \cite{JK95}, Jeffrey and Kirwan proved a formula (Theorem 8.1 in \cite{JK95}) computing the following
cohomological quantity
\begin{equation}
\kappa_{0}(\eta) e^{\I \omega_{\red}} [M_{\red}]
\end{equation}
for any $ \eta \in H^{*}_{G}(M) $. Also in \cite{JK3}, they rewrote the residue formula
(Theorem 3.1 in \cite{JK3}) computing the following cohomological quantity
\begin{equation}
\kappa_{0}(\eta) e^{\omega_{\red}} [M_{\red}]
\end{equation}
for any $ \eta \in H^{*}_{G}(M) $.
Both articles basically compute the same quantity.

Before we state the formula, we first fix some notations which will be convenient for
later use.

The following notations can be applied to any general semisimple
compact connected Lie group $ G $. In this section, to avoid
ambiguity, we assume $ G = \SU(3) $ throughout.
\begin{itemize}
\item Let $ s $ denote the real dimension of $ G $.
\item Let $ l $ denote the real dimension of $ T $.
\item Let $ \varpi $ denote the product of the positive roots of $ G $ (here, roots are regarded as elements in $ \lietdual $), that is,
\begin{equation}
\varpi(X) = \prod_{\gamma \in \proots} \gamma(X)
\end{equation}
for all $ X $ in $ \liet $ or $ \liet_{\CC} $, where $ \liet_{\CC} $ denotes
the complexification of $ \liet $. Therefore we can regard $ \varpi $ as
a polynomial function on $ \liet $ or $ \liet_{\CC} $. 
\item Notice that
\begin{equation}
\varpi(w \cdot X) = \sgn(w) \varpi(X)
\end{equation}
for all $ w \in W $ and all $ X $ in $ \liet $ or $ \liet_{\CC} $.
\item The fixed $ G $-invariant inner product $ (\cdot, \cdot) $ induces measures on $ \lieg $ and $ \liet $
and their corresponding complexifications, $ \lieg_{\CC} $ and $ \liet_{\CC} $. 
\item Following \cite{JK95},
we will reserve the Greek letter $ \phi $ for a variable in $ \lieg $ or $ \lieg_{\CC} $ and the Greek letter
$ \psi $ for a variable in $ \liet $ or $ \liet_{\CC} $.
\item Let $ [d \phi] $ denote the induced measure on $ \lieg $ or $ \lieg_{\CC} $.
\item Let $ [d \psi] $ denote the induced measure on $ \liet $ or $ \liet_{\CC} $.
\item The inner product also induces Riemannian volume forms on $ G $ and $ T $.
\item Let $ \vol^{\RRR}(G) $ denote the Riemannian volume of $ G $ under this induced Riemannian volume form on
$ G $.
\item Let $ \vol^{\RRR}(T) $ denote the Riemannian volume of $ T $ under this induced Riemannian volume form on
$ T $.
\item Let $ \mu_{T}: M \to \liet $ denote the composition of the moment map $ \mu_{G} $ with the orthogonal
projection from $ \lieg $ to $ \liet $. Notice that the orthogonal projection from $ \lieg $ to $ \liet $
corresponds to the restriction map from $ \liegdual $ to $ \lietdual $ under the identification through
the inner product. Therefore, $ \mu_{T} $ is a moment map for the action of $ T $ on $ M $.
\end{itemize}

Now we can state the residue formula of Jeffrey and Kirwan, adapted to our situation. 
\begin{theorem}{\cite{JK95} Theorem 8.1, \cite{JK3} Theorem 3.1} \label{thm:residue}
Let $ G $ be a general semisimple
compact connected Lie group. Let $ M = \orb_{\xi_{1}} \times \dots \times \orb_{\xi_{N}} $
be the Cartesian product of $ N \geq 3 $ adjoint orbits of $ G $, where the $ \xi_{i} $
satisfy the assumptions (A1), (A2) and (A3) outlined in \S \ref{sss:assumptions}. Then,
for all $ \eta \in H^{*}_{G}(M) $, we have
\begin{equation}\label{eq:residue}
\kappa_{0}(\eta) e^{\I \omega_{\red}} [M_{\red}]
= n_{0} C_{G} \res\left( \varpi^{2}(\psi) \sum_{F \in M^{T}} r^{\eta}_{F}(\psi) [d \psi] \right) 
\end{equation}
where $ n_{0} $ is the number of points in  the stabilizer subgroup in $ G $ of a generic point in $ M_{0} $,
and the constant $ C_{G} $ is defined by
\begin{equation}
C_{G} := \frac{(-1)^{n_{+}}}{(2 \pi)^{s - l} \left| W \right| \vol^{\RRR}(T)}.
\end{equation}
Here, $ n_{+} $ denotes the number of positive 
roots, that is, $ n_{+} = (s - l) / 2 $.
$ \psi $ is a variable in $ \liet_{\CC} $. Notice that in our situation, all fixed points in $ M^{T} $
are isolated. If $ F \in M^{T} $, the meromorphic function $ r^{\eta}_{F} $ on $ \liet_{\CC} $ is
defined by
\begin{equation}
r^{\eta}_{F}(\psi) := e^{\I (\mu_{T}(F), \psi)} \int_{F} \frac{i_{F}^{*}(\eta(\psi) e^{\I \omega})}{\ee_{F}(\psi)}
\end{equation}
where $ i_{F}: F \to M $ is the inclusion and $ \ee_{F} $ is the $ T $-equivariant Euler class of the
normal bundle to $ F $ in $ M $. The notation ${\rm res}$ will be defined 
in Section \ref{ss:residue} below.
\end{theorem}

\begin{remark}
Recall that for each point $ \veceta \in M_{0} $, $ \Stab_{G}(\veceta) $ is finite.
A generic point in $ M_{0} $ is a point $ \veceta \in M_{0} $ such that
the cardinality of $ \Stab_{G}(\veceta) $ is the smallest among all points in $ M_{0} $. Therefore,
the number $ n_{0} $ above is well defined. In particular, the above theorem
does not require the action $ G \circlearrowright M_{0} $ to be free.
\end{remark}

Since $ F $ is an isolated fixed point, we have
\begin{equation}
\int_{F} \frac{i_{F}^{*}(\eta(\psi) e^{\I \omega})}{\ee_{F}(\psi)}
= \frac{i_{F}^{*}(\eta(\psi) e^{\I \omega})}{\ee_{F}(\psi)}
\end{equation}
and the quantity
\begin{equation}
i_{F}^{*}(\eta(\psi) e^{\I \omega})
\end{equation}
will be a $ \CC $-valued function on $ F $, that is, a complex number. We denote this complex number by
\begin{equation}
\eta_{F}(\psi) := i_{F}^{*}(\eta(\psi) e^{\I \omega}) \in \CC.
\end{equation}

Furthermore, since $ M = \orb_{\xi_{1}} \times \orb_{\xi_{2}} \times \orb_{\xi_{3}} $,
any $ F \in M^{T} $ can be written (see Proposition \ref{prop:TfixedPoints.N=3}) by
\begin{equation}
F = (w_{1} \cdot \xi_{1}, w_{2} \cdot \xi_{2}, w_{3} \cdot \xi_{3})
\end{equation}
for some $ \vecw = (w_{1}, w_{2}, w_{3}) \in W \times W \times W $.
Then we can write
\begin{equation}
F = \vecw \cdot \vecxi := (w_{1} \cdot \xi_{1}, w_{2} \cdot \xi_{2}, w_{3} \cdot \xi_{3}).
\end{equation}
In this case, $ \ee_{F}(\psi) $ is the complex number
\begin{equation}
\sgn(\vecw) \varpi^{3}(\psi)
\end{equation}
where
\begin{equation}
\sgn(\vecw) := \prod_{i = 1}^{3} \sgn(w_{i}).
\end{equation}
Recall that $\varpi(\psi) = \prod_{\gamma > 0} \gamma(\psi)$ is the 
product of all positive roots  evaluated on $\psi \in \liet$.
In this case, we can define
\begin{equation}
\sgn(F) := \sgn(\vecw).
\end{equation}
Here $\sgn(w)$  is the signature of the Weyl group element $w$  -- this is equal to its signature
as a permutation.
In addition, we write $ \vecw(F) $ to mean the $ \vecw $ such that $ F = \vecw \cdot \vecxi $.

Now, we can rewrite Equation (\ref{eq:residue}) as
\begin{equation}
\kappa_{0}(\eta) e^{\I \omega_{\red}} [M_{\red}]
= n_{0} C_{G} \res\left( \varpi^{2}(\psi) \sum_{F \in M^{T}} e^{\I (\mu_{T}(F), \psi)}  \frac{\eta_{F}(\psi)}{\sgn(F) \varpi^{3}(\psi)} [d \psi] \right) .
\end{equation}
This can be further simplified to
\begin{equation}
\kappa_{0}(\eta) e^{\I \omega_{\red}} [M_{\red}]
= n_{0} C_{G} \res\left( \sum_{F \in M^{T}} \sgn(F) \frac{\eta_{F}(\psi) e^{\I (\mu_{T}(F), \psi)}}{\varpi(\psi)} [d \psi] \right) .
\end{equation}

We have the following result as a corollary of Theorem \ref{thm:residue}.

\begin{corollary}\label{cor:vol}
Let $ M = \orb_{\xi_{1}} \times \orb_{\xi_{2}} \times \orb_{\xi_{3}} $ be the Cartesian
product of $ N = 3 $ adjoint orbits of $ G = \SU(3) $, where the $ \xi_{i} $
 satisfy
the assumptions (A1), (A2) and (A3) outlined in \S \ref{sss:assumptions}. Then,
the symplectic volume of $ M_{\red} $ can be expressed as
\begin{equation}
\vol^{\SSS}(M_{\red}) = \frac{1}{\I^{d / 2}} n_{0} C_{G} \res\left( \sum_{\vecw \in W^{3}} \sgn(\vecw) \frac{e^{\I (w_{1} \cdot \xi_{1} + w_{2} \cdot \xi_{2} + w_{3} \cdot \xi_{3}, \psi)}}{\varpi(\psi)} [d \psi] \right) .
\end{equation}
Here, $ d $ is the real dimension of $ M_{\red} $. 

\end{corollary}

\begin{proof}
To obtain the symplectic volume $ \vol^{\SSS}(M_{\red}) $, we let $ \eta = 1 \in H^{*}_{G}(M) $
and compute:
\begin{align}
\vol^{\SSS}(M_{\red}) &= \frac{1}{\I^{d / 2}} e^{\I \omega_{\red}} [M_{\red}]  \notag \\
&= \frac{1}{\I^{d / 2}} \kappa_{0}(1) e^{\I \omega_{\red}} [M_{\red}] \notag \\
&= \frac{1}{\I^{d / 2}} n_{0} C_{G} \res\left( \sum_{F \in M^{T}} \sgn(F) \frac{e^{\I (\mu_{T}(F), \psi)}}{\varpi(\psi)} [d \psi] \right) .
\end{align}

Since the fixed point set $ M^{T} $ can be parametrized by
\begin{equation}
\vecw = (w_{1}, w_{2}, w_{3}) \in W^{3},
\end{equation}
the sum
\begin{equation}
\sum_{F \in M^{T}} \sgn(F) \frac{e^{\I (\mu_{T}(F), \psi)}}{\varpi(\psi)}
\end{equation}
can be rewritten as
\begin{equation}
\sum_{\vecw \in W^{3}} \sgn(\vecw) \frac{e^{\I (w_{1} \cdot \xi_{1} + w_{2} \cdot \xi_{2} + w_{3} \cdot \xi_{3}, \psi)}}{\varpi(\psi)}.
\end{equation}
\end{proof}

To fully understand the above formulas, 
we need to understand the residue  $ \res $. We describe this
in the next section.

\section{The Residue} \label{s:residue}

In this section  we shall study the residue  $ \res $ in detail. The main references for this
section are \cite{JK95} and \cite{JK3}.

\subsection{Definition of the residue} \label{ss:residue}
The Fourier transform is important to the development of  nonabelian localization in \cite{JK95}.
Following \cite{JK95}, we reserve the letter $ z $ for a variable in $ \liegdual $ and the letter
$ y $ for a variable in $ \lietdual $.

If $ f: \lieg \to \CC $ is a tempered distribution, we define its Fourier transform $ F_{G}f $ on $ \liegdual $ to be
\begin{equation}
(F_{G}f)(z) := \frac{1}{(2 \pi)^{s / 2}} \int_{\phi \in \lieg} f(\phi) e^{-\I \left\langle z, \phi \right\rangle } [d\phi].
\end{equation}
Since $ f: \lieg \to \CC $ can be regarded as a tempered distribution on $ \liet $ by restriction
to $ \liet $, we can also define its Fourier transform $ F_{T}f $ on $ \lietdual $ to be
\begin{equation}
(F_{T}f)(y) := \frac{1}{(2 \pi)^{l / 2}} \int_{\psi \in \liet} f(\psi) e^{-\I \left\langle y, \psi \right\rangle } [d\psi].
\end{equation}


Proposition 8.4 in \cite{JK95} 
(see also \cite{Hormander1})\label{prop:8.4} \hfill
gives a characterization of the residue.
Based on this proposition, we have the following definition:

\begin{definition}[Definitions 8.5, the paragraph before Definition 8.8, and Definition 8.8
in \cite{JK95}]\label{def:8.5,8.9}
Let $ \Lambda $ be a proper cone in $ \liet $. Let $ \Lambda^{0} $ denote the interior of
 $ \Lambda $. Let $ h $ be a holomorphic function on
$ \liet - \I \Lambda^{0} \subset \liet_{\CC} $ such that for every compact subset $ K $ 
of $ \liet - \I \Lambda^{0} $
there exists a constant $ C_{K} $ and an integer $ N_{K} \geq 0 $ such that
\begin{equation}
\left| h(\zeta) \right| \leq C_{K} (1 + \left| \zeta \right| )^{N_{K}}
\end{equation}
for all $ \zeta \in K $. Let $ \chi: \liegdual \to \RR $ be a smooth invariant function 
with
compact support and strictly positive in some neighbourhood of $ 0 $. Let $ \hat{\chi} = 
F_{G}\chi:
\lieg \to \CC $ be its Fourier transform. Let $ \hat{\chi}_{\epsilon}(\phi) = \hat{\chi}(\epsilon \phi) $
so that
\begin{equation}
(F_{G}\hat{\chi}_{\epsilon})(z) = \chi_{\epsilon}(z) := \frac{1}{\epsilon^{s}} \chi(\frac{z}{\epsilon}).
\end{equation}
Also we assume $ \hat{\chi}(0) = 1 $. Then we define
\begin{equation}\label{def:res}
\res^{\Lambda, \chi}(h(\psi) [d\psi]) := \lim_{\epsilon \to 0^{+}} \frac{1}{(2 \pi \I)^{l}}
\int_{\psi \in \liet - \I \xi} \hat{\chi}(\epsilon \psi) h(\psi) [d\psi],
\end{equation}
where $ \xi $ is any element of $ \Lambda^{0} $. Furthermore, in the case where $ h $
is a sum of other functions $ h_{i} $ and $ F_{T}h $ is smooth at $ 0 $ but the Fourier 
transforms of
$ h_{i} $ may not be smooth at $ 0 $, we need to introduce a small generic parameter 
$\rho \in \lietdual $
so that all the functions in this sum have Fourier transforms that are smooth at $ 0 $. 
More precisely,
let $ \Lambda $, $ \chi $ and $ h $ be as in the above. Let $ \rho \in \lietdual $ be 
such that the
distribution $ F_{T}h $ is smooth on the ray $ t \rho $ for $ t \in (0, \delta) $ for 
some $ \delta > 0 $,
and suppose $ (F_{T}h)(t \rho) $ tends to a well defined limit as $ t \to 0^{+} $. Then
we define
\begin{equation}
\res^{\rho, \Lambda, \chi}(h(\psi) [d\psi]) := \lim_{t \to 0^{+}} \res^{\Lambda, \chi}(
h(\psi)
e^{\I \left\langle t \rho, \psi \right\rangle } [d\psi]).
\end{equation}
\end{definition}

\begin{remark}
By the paragraph after Definition 8.5 in \cite{JK95}, the integral in (\ref{def:res}) converges and
is independent of $ \xi \in \Lambda^{0} $. Furthermore, by Propositions 8.6, 8.7 and 8.9 \
in \cite{JK95},
the residue of a meromorphic form $ \Omega $ is independent of the choices of $ \Lambda $\
, $ \chi $ and $ \rho $
if $ \Omega $ is sufficiently well behaved, as is the case for Theorem \ref{thm:residue}.
\end{remark}

\subsection{Nonabelian localization}

Before we go into the computational aspects of the residue , we shall briefly
sketch the proof of Theorem \ref{thm:residue} by summarizing key points in \cite{JK95}, so that
we will have a better idea about why the nonabelian localization technique in \cite{JK95} works.

The first key point is the abelian localization formula:

\begin{theorem}[Berline and Vergne, \cite{AB}, \cite{Audin}, \cite{BV83}, \cite{DH}; Theorem 2.1 in \cite{JK95}]\label{thm:ab}
Let $ G $ be a general semisimple
compact connected Lie group and $ M $
a symplectic manifold equipped with a Hamiltonian action of $G$. 
Assume that $0$ is a regular value of the moment map. 
If $ \sigma \in H^{*}_{T}(M) $, then
\begin{equation}
(\Pi_{*} \sigma)(\psi) = \sum_{F \in M^{T}} \int_{F} \frac{i_{F}^{*}(\sigma(\psi))}{\ee_{F}(\psi)}.
\end{equation}
\end{theorem}

A special case of Theorem \ref{thm:ab} is the case when $M$ is a product of coadjoint orbits:
\begin{theorem}\label{thm:absp}
Let $ G $ be a general semisimple
compact connected Lie group and $ M = \orb_{\xi_{1}} \times \dots \times \orb_{\xi_{N}} $
be the Cartesian product of $ N \geq 3 $ adjoint orbits of $ G $, where the $ \xi_{i} $
satisfy the assumptions (A1), (A2) and (A3) outlined in \S \ref{sss:assumptions}.
If $ \sigma \in H^{*}_{T}(M) $, then
\begin{equation}
(\Pi_{*} \sigma)(\psi) = \sum_{F \in M^{T}} \int_{F} \frac{i_{F}^{*}(\sigma(\psi))}{\ee_{F}(\psi)}.
\end{equation}
\end{theorem}

We are interested in the case when $ \sigma = \eta e^{\I \bar{\omega}} $ where $ \eta \in H_{G}^{*}(M) $ and
\begin{equation}
\bar{\omega}(\phi) = \omega + \mu_{G}(\phi)
\end{equation}
is the standard equivariant extension of $ \omega $.
Notice that here $ \sigma = \eta e^{\I \bar{\omega}} $ is regarded as a $ T $-equivariant cohomology class
through the restriction map $ \liegdual \to \lietdual $.

Let $ r^{\eta} := \Pi_{*}(\eta e^{\I \bar{\omega}}) \in H^{*}_{T} $,
where we have defined $H^*_T:= H^*_T ({\mathrm pt})$ 
as the $T$-equivariant cohomology of a point. This was (\ref{ss:eqc}). By Theorem \ref{thm:ab}, we have
\begin{equation}\label{eq:r.eta}
r^{\eta}(\psi) = \sum_{F \in M^{T}} r^{\eta}_{F}(\psi),
\end{equation}
where
\begin{equation}\label{eq:r.eta.f}
r^{\eta}_{F}(\psi) = e^{\I (\mu_{T}(F), \psi)} \int_{F} \frac{i_{F}^{*}(\eta(\psi) e^{\I \omega})}{\ee_{F}(\psi)}.
\end{equation}

It turns out that the Fourier transform of $ \Pi_{*}(\eta e^{\I \bar{\omega}}) $ will be
closely related to the cohomological quantity $ \kappa_{0}(\eta) e^{\I \omega_{\red}} [M_{\red}] $.

Proposition 8.10 in \cite{JK95} gives an explicit characterization of the residue.

\begin{remark}

The rest of the proof for \cite{JK95} Proposition 8.10 depends on the following important result, the normal form theorem for the symplectic structure and
moment map near a regular value. The normal form theorem is due to 
 Gotay \cite{Gotay}, Guillemin and Sternberg \cite{GS} and Marle \cite{Marle}. 
The statement is as follows.
\begin{proposition}[Proposition 5.2 in \cite{JK95}]\label{prop:nf}
Assume $ 0 $ is a regular value of $ \mu_{G} $. Then there is a neighbourhood
\begin{equation}
U \cong M_{0} \times \left\lbrace z \in \liegdual \  : \  \left| z \right| < h \right\rbrace \subset
M_{0} \times \liegdual,
\end{equation}
where $ h > 0 $ is some sufficiently small number, of $ M_{0} $ on which the symplectic form $ \omega $ can be given as follows.
Recall that $ p_{0}: M_{0} \to M_{\red} $ is the orbifold principal $ G $-bundle.
Let $ \theta \in \Omega^{1}(M_{0}) \otimes \lieg $ be a connection on this principal bundle.
Recall that on $ M_{\red} $ there is a symplectic structure $ \omega_{\red} $ such that
$ p_{0}^{*} \omega_{\red} = i_{0}^{*} \omega $. Let $ \alpha $ be a $ 1 $-form on $ U \subset M_{0} \times
\liegdual $ defined by
\begin{equation}
\alpha_{(p, z)}(v, \xi) := \left\langle z, \theta_{p}(v) \right\rangle 
\end{equation}
for all $ p \in M_{0} $, $ z \in \liegdual $ with $ \left| z \right| < h $, $ v \in T_{p}M_{0} $ and $ \xi \in T_{z}\liegdual = \liegdual $.
Then the symplectic form $ \omega $ on $ U $ is 
\begin{equation}
\omega = \pr_{1}^{*} p_{0}^{*} \omega_{\red} + d \alpha.
\end{equation}
Moreover, the moment map $ \mu_{G} $ on $ U $ is  $ \mu_{G}(p, z) = z $.
\end{proposition}
\end{remark}

The authors of \cite{JK95} used  a sequence of appropriately chosen test functions $ \chi_{\epsilon}: \liegdual \to \RR_{\geq 0} $ such
that as $ \epsilon \to 0 $, the functions $ \chi_{\epsilon} $ tend to the Dirac delta distribution 
on $ \liegdual $.  They  integrated $ F_{G}(\Pi_{*}(\eta e^{\I \bar{\omega}})) $ against this sequence
of test functions $ \chi_{\epsilon} $. By invoking Proposition \ref{prop:nf}, they were able to concentrate
on arbitrarily small neighbourhoods of $ M_{0} $ and obtain the estimate which eventually established the
link between $ F_{G}(\Pi_{*}(\eta e^{\I \bar{\omega}}))(0) $ and the cohomological quantity
$ \kappa_{0}(\eta) e^{\I \omega_{\red}} [M_{\red}] $, finishing the proof for \cite{JK95} Proposition 8.10.



We need two more properties of the residue. 
The first is as follows.
\begin{proposition}[Proposition 8.7 in \cite{JK95}]\label{prop:8.7}
Let $ u: \lietdual \to \CC $ be a distribution, and assume the set $ \Gamma_{u} $ defined in
Proposition \ref{prop:8.4} contains $ -\Lambda^{0} $. Then, $ h = F_{T}u $ is a holomorphic function
on $ \liet - \I \Lambda^{0} $ and $ h $ satisfies the hypotheses in Definition \ref{def:8.5,8.9}.
Assume in addition that $ u $ is smooth at $ 0 $. Then $ \res^{\Lambda, \chi}(h(\psi)[d\psi]) $ is
independent of the test function $ \chi $, and moreover,
\begin{equation}
\res^{\Lambda, \chi}(h(\psi)[d\psi]) = \frac{1}{\I^{l} (2\pi)^{l/2}} u(0).
\end{equation}
\end{proposition}

 The second important realization is provided by the following
point due to Guillemin, Lerman, Prato and Sternberg (\cite{GLS1}, \cite{GLS} and \cite{GP}):

\begin{proposition}[Proposition 3.6 in \cite{JK95}]\label{prop:3.6} \hfill
\begin{description}
\item[(a)] (Part (a) of Proposition 3.6 in \cite{JK95}; see \cite{GLS1}, Section 3.2 in \cite{GLS}, and \cite{GP}) Define
\begin{equation}
H_{\bar{\beta}}(y) = \vol\left\lbrace (s_{1}, \dots, s_{\nu}) \  : \  s_{i} \geq 0, \  y = \sum_{j} s_{j}
\beta_{j} \right\rbrace 
\end{equation}
for some $ \nu $-tuple $ \bar{\beta} = (\beta_{1}, \dots, \beta_{\nu}) $
with $ \beta_{j} \in \lietdual $ such that all $ \beta_{j} $ lie in the interior of some half-space of $ \lietdual $, where $\nu$ is a positive integer.

 Here the reason for $ s_{i} \geq 0 $ is that we are working with the cone
spanned by $ \bar{\beta} $.
Here $ \vol $ denotes the standard Euclidean volume multiplied by a normalization constant
which is chosen so that Equation (\ref{eq:FT.H}) below holds.
Thus, $ H_{\bar{\beta}} $ is a piecewise polynomial function supported on the cone
\begin{equation}
C_{\bar{\beta}} := \left\lbrace \sum_{j} s_{j} \beta_{j} \  : \  s_{j} \geq 0 \right\rbrace .
\end{equation}
Let $ h(y) := H_{\bar{\beta}}(y + \tau) $ for some $ \tau \in \lietdual $.
Then the Fourier transform of $ h $ is given for $ \psi $ in the complement of the union of the hyperplanes
$ \left\lbrace \psi \in \liet \  : \  \beta_{j}(\psi) = 0 \right\rbrace  $ by the formula
\begin{equation}\label{eq:FT.H}
F_{T}h(\psi) = \frac{e^{\I \left\langle \tau, \psi \right\rangle }}{\I^{\nu} \prod_{j = 1}^{\nu}\beta_{j}(\psi)}.
\end{equation}
\end{description}
\end{proposition}

By these considerations, $ \res(r^{\eta}_{F}(\psi) [d\psi]) $ can be computed as
\begin{equation}
\res^{\rho, \Lambda, \chi}(r^{\eta}_{F}(\psi) [d\psi]) = \lim_{t \to 0^{+}}
\frac{1}{(2\pi)^{l/2} \I^{l}} F_{T}r^{\eta}_{F}(t \rho).
\end{equation}
Recall that $r^\eta$ was defined above in (\ref{eq:r.eta}),
by Definition \ref{def:8.5,8.9} and Proposition \ref{prop:8.7}. Thus
\begin{align}
\sum_{F \in M^{T}} \res^{\rho, \Lambda, \chi}(\varpi^{2}(\psi) r^{\eta}_{F}(\psi) [d\psi]) &=
\res^{\rho, \Lambda, \chi}(\varpi^{2}(\psi) r^{\eta}(\psi) [d\psi]) \\
&= \frac{1}{(2\pi)^{l/2} \I^{l}} F_{T}(\varpi^{2} r^{\eta})(0),
\end{align}
providing the last link for the proof of Theorem \ref{thm:residue}.

Please note the following.

\begin{enumerate}
\item  Theorem 8.1, \cite{JK95}:  $$i \vol(M)  = e^{i \omega_0}[\mathcal{M}_X]   = (-1)^{n_+ }B  Res  \sum_w (-1)^w e^{i <w 
\lambda, X>} /\mathcal{D}(X) $$
where $B$ is a real constant.
The constant $B$ equals
$$ \frac{1}{(2 \pi)^{s-\ell} |W| \vol(T)}. $$
\item  \cite{JK95} Proposition 8.11 (ii): $ 
Res (\frac{e^{i \lambda(X)}}{\prod_j \beta_j(X)}   = i^{n_+-\ell}  H_\beta$
where $  H_\beta $ was defined in \cite{JK95} Proposition 3.6 
 equations (3.18). See also \cite{JK95} (8.28),
where it is stated that the residue of a function $h$ is equal to 
$$\frac{1}{2 \pi i } \int h(\psi) d\psi.$$
In the $SU(3)$ case we have ${n_+ - \ell}  = 3-2 =1$. 
In other words, in ourcase
 there is a multiplicative factor of $i$ in the definition of
the residue.
\end{enumerate}

The factor $i$ on the right hand side comes from $\sum_w (-1)^w e^{i <w \lambda, X>}.$
Using the Weyl character formula this is proportional to $i^{n_+}$ times a real number (where $n_+$ 
is the  number of positive
roots).
In our case $n_+  = 3$.
So when $G  = SU(3)$, we  
get $i {\rm vol} (M) = i^{n_+} $ times a real number. So we have identified the volume.

\subsection{Computing the residue}

In this section, we will focus on computing the residue $ \res(\Omega_{\lambda}) $ of a special class of
meromorphic forms $ \Omega_{\lambda} $ such that
\begin{equation}
\Omega_{\lambda}(\psi) = \frac{e^{\I \lambda(\psi)} [d\psi]}{\prod_{j = 1}^{\nu} \beta_{j}(\psi)},
\end{equation}
where $ \lambda $ is some point in $ \lietdual $ and $ \beta_{j} $ all lie in the dual cone of a proper cone $ \Lambda $ in $ \liet $. A proper cone is
defined as  an open cone such that its apex is the origin and it is properly
contained in some half space. Given a proper cone $ \Lambda \subset \liet $, its dual cone $ \Lambda^{*} $ is defined to be the
following collection of elements in $ \lietdual $:
\begin{equation}
\Lambda^{*} := \left\lbrace \beta \in \lietdual \  : \  \beta(\psi) > 0 \  \text{for all} \  \psi \in \Lambda \right\rbrace .
\end{equation}

In Proposition 8.11 in \cite{JK95},
a list of properties satisfied by the residue is given
  (see also Proposition 3.2 in \cite{JK3}). This list of properties will be the basis
for the computations in this article.

Now, we are finally ready for the computation of the symplectic volume of the triple ($ N = 3 $) reduced product $ M_{\red} $, i.e., the cohomological quantity
\begin{equation}
\frac{1}{\I^{d/2}} e^{\I \omega_{\red}} [M_{\red}].
\end{equation}

We have the following result.

\begin{theorem}\label{thm:trp.vol} 
The symplectic volume of the triple reduced product $ M_{\red}(\vecxi) $ of $ G = \SU(3) $ (with
the quantities $ \vecxi = (\xi_{1}, \xi_{2}, \xi_{3}) $ satisfying (A1), (A2) and (A3) outlined in \S \ref{sss:assumptions}) 
can be computed by the following explicit formula:
\begin{equation}\label{eq:trp.vol}
\vol^{\SSS}(M_{\red}(\vecxi)) = K \cdot \sum_{i = 0}^{5} \sum_{j = 0}^{5} \sum_{k = 0}^{5}
(-1)^{i+j+k} \cdot \max\left( \min\left( P_{ijk}^{(1)}(\vecxi), P_{ijk}^{(2)}(\vecxi) \right), 0 \right) ,
\end{equation}
where $ K $ is a real  constant that depends on the inner product and $ \bar{\beta} $
(although the overall formula does not depends on these choices) and
\begin{align}
P_{ijk}^{(1)}(\vecxi) &= \left( \frac{2}{3} \pr_{1} - \frac{1}{3} \pr_{2}\right) \left( P_{ijk}(\vecxi)\right) , \\
P_{ijk}^{(2)}(\vecxi) &= \left( \frac{1}{3} \pr_{1} + \frac{1}{3} \pr_{2}\right) \left( P_{ijk}(\vecxi)\right) ,
\end{align}
with
\begin{equation}
P_{ijk}(\vecxi) = 
\ws_{i} \cdot
\begin{pmatrix}
\ell_{1} \\
m_{1}
\end{pmatrix}
+
\ws_{j} \cdot
\begin{pmatrix}
\ell_{2} \\
m_{2}
\end{pmatrix}
+
\ws_{k} \cdot
\begin{pmatrix}
\ell_{3} \\
m_{3}
\end{pmatrix} \in \RR^{2}.
\end{equation}
Here $ \pr_{i} $ denotes the standard projection onto the $ i $-th coordinate.
 Here
\begin{equation}
\xi_{i} = (\ell_{i} - m_{i}) \cdot \Omega_{1} + m_{i} \cdot \Omega_{2} = \ell_{i} \cdot \Omega_{1} + 
m_{i} \cdot (\Omega_{2} - \Omega_{1}), \  \ell_{i} > m_{i} > 0,
\end{equation}
\end{theorem}

\begin{proof}
By Corollary \ref{cor:vol}, we have
\begin{equation}\label{eq:vol.res}
\frac{1}{\I^{d/2}} e^{\I \omega_{\red}} [M_{\red}]
= \frac{1}{\I^{d / 2}} n_{0} C_{G} \res\left( \sum_{\vecw \in W^{3}} \sgn(\vecw) \frac{e^{\I (w_{1} \cdot \xi_{1} + w_{2} \cdot \xi_{2} + w_{3} \cdot \xi_{3}, \psi)}}{\varpi(\psi)} [d \psi] \right) .
\end{equation}

To compute the residue on the right hand side of the above equation,
we need to first choose an open cone $ \Lambda $ in $ \liet $.
We choose $ \Lambda = \opwc $, the open positive Weyl chamber.
One reason that we make this choice is that we observe that (here recall that we are considering $ G = \SU(3) $)
\begin{equation}
\varpi(\psi) = \prod_{j = 1}^{3} \beta_{j}(\psi),
\end{equation}
where
\begin{align}
\beta_{1} &= (\diag(2 \pi \I, -2 \pi \I, 0), \cdot), \\
\beta_{2} &= (\diag(0, 2 \pi \I, -2 \pi \I), \cdot), \\
\beta_{3} &= (\diag(2 \pi \I, 0, -2 \pi \I), \cdot).
\end{align}
Thus $ \beta_{3} = \beta_{1} + \beta_{2} $ and the collection $ \left\lbrace \beta_{1}, \beta_{2}, \beta_{3} \right\rbrace  $ is just the set of positive roots of $ G = \SU(3) $.
Notice that all of $ \beta_{1}, \beta_{2}, \beta_{3} $ lie in $ \Lambda^{*} $, the dual cone of $ \Lambda $.

Let $ \bar{\beta} = (\beta_{1}, \beta_{2}, \beta_{3}) $.

Let
\begin{equation} \label{circledef}
\vecw \odot \vecxi := \sum_{i = 1}^{N} w_{i} \cdot \xi_{i}.
\end{equation}
Here we are considering $ N = 3 $. 

Let us rewrite the residue part in Equation (\ref{eq:vol.res}):
\begin{align}
&\res\left( \sum_{\vecw \in W^{3}} \sgn(\vecw) \frac{e^{\I (w_{1} \cdot \xi_{1} + w_{2} \cdot \xi_{2} + w_{3} \cdot \xi_{3}, \psi)}}{\varpi(\psi)} [d \psi] \right) \label{eq:res.totalsum.N=3} \\
= &\res \left( \sum_{\vecw \in W^{3}} \sgn(\vecw) \frac{e^{\I (\vecw \odot \vecxi)(\psi)} [d\psi]}{\prod_{j = 1}^{3} \beta_{j}(\psi)} \right) \\
= &\sum_{\vecw \in W^{3}} \sgn(\vecw) \res \left( \frac{e^{\I (\vecw \odot \vecxi)(\psi)} [d\psi]}{\prod_{j = 1}^{3} \beta_{j}(\psi)} \right) \\
= &\sum_{\vecw \in W^{3}} \sgn(\vecw) \frac{\I^{3}}{(2 \pi \I)^{2}} H_{\bar{\beta}}(\vecw \odot \vecxi).
\end{align}

Combining the above formula with the constant part in Equation (\ref{eq:vol.res}) and
recalling that 
$n_+ = 3$, we have
$$C_G = -\frac{1}{(2\pi)^4 \cdot 6  \cdot \vol(T)} $$
so 

\begin{equation} \label{eq:trp.res3} 
 \vol^{\SSS}(M_{\red}(\vecxi)) = \frac{1}{\I^{d/2}} n_{0} C_{G} \frac{\I^{3}}{(2 \pi \I)^{2}}
\sum_{\vecw \in W^{3}} \sgn(\vecw) H_{\bar{\beta}}(\vecw \odot \vecxi)  \end{equation}
$$=   -3 \cdot \frac{1}{(2 \pi)^{6}   \cdot 6 \cdot \vol^{\RRR}(T) } \cdot
\frac{1}{(2 \pi )^{2} } \cdot \sum_{\vecw \in W^{3} } \sgn(\vecw) H_{\bar{\beta}}(\vecw \odot \vecxi)  $$
$$= -\frac{1}{2 \cdot (2 \pi)^{8} \cdot \vol^{\RRR}(T)} \sum_{\vecw \in W^{3}} \sgn(\vecw) H_{\bar{\beta}}(\vecw \odot \vecxi) . $$

Notice that since we are considering triple reduced products of $ G = \SU(3) $ here, we have
\begin{align}
d &= N (s - l) - 2s = 3 \cdot (8 - 2) - 2 \cdot 8 = 2 ,\\
n_{+} &= (s - l) / 2 = (8 - 2) / 2 = 3, \\
n_{0} &= \left| \ZZ(\SU(3)) \right| = 3.
\end{align}
Note that $d$ is the real dimension of the reduced space. 
Therefore we need to compute $ H_{\bar{\beta}}(\vecw \odot \vecxi) $ to obtain an explicit formula
for the symplectic volume of a triple reduced product.
Since later we will compare our volume formula with the volume formula obtained by
Suzuki and Takakura \cite{ST08}, we will parametrize things in a way similar to theirs.


Any element $ \xi \in \opwc $ can be written as
\begin{equation}
\xi = (\ell - m) \cdot \frac{2 \beta_{1} + \beta_{2}}{3} + m \cdot \frac{\beta_{1} + 2 \beta_{2}}{3}
\end{equation}
for some $ \ell > m > 0 $. If $ \ell $ and $ m $ can be any real number, then the above formula
parametrizes all $ \xi \in \liet $.

Let
\begin{align}
\Omega_{1} &:= \frac{2 \beta_{1} + \beta_{2}}{3}, \\
\Omega_{2} &:= \frac{\beta_{1} + 2 \beta_{2}}{3}.
\end{align}
Then
\begin{equation}
\xi = (\ell - m) \cdot \Omega_{1} + m \cdot \Omega_{2} = \ell \cdot \Omega_{1} + m \cdot (\Omega_{2} - \Omega_{1}).
\end{equation}

We fix $ \left\lbrace \Omega_{1}, \Omega_{2} - \Omega_{1} \right\rbrace  $ as the basis for $ \liet $
for this computation. By using this basis, it will be easier for us to compare our volume formula
with the volume formula obtained by Suzuki and Takakura.
Our goal here is to express $ H_{\bar{\beta}}(\xi) $ as a function of $ \ell $ and $ m $.
It is helpful to express $ H_{\bar{\beta}}(\lambda_{1} \beta_{1} + \lambda_{2} \beta_{2}) $
in terms of $ \lambda_{1} $ and $ \lambda_{2} $ first.

By the definition of $ H_{\bar{\beta}} $, we have
\begin{equation}
H_{\bar{\beta}}(\lambda_{1} \beta_{1} + \lambda_{2} \beta_{2}) = 
\vol \left\lbrace (s_{1}, s_{2}, s_{3}) \in \RR_{+}^{3} \  : \  \sum_{j = 1}^{3} s_{j} \beta_{j}
= \lambda_{1} \beta_{1} + \lambda_{2} \beta_{2} \right\rbrace ,
\end{equation}
where $ \RR_{+} $ denotes the set of all nonnegative real numbers.
Therefore we want to solve the following equation:
\begin{equation}
s_{1} \beta_{1} + s_{2} \beta_{2} + s_{3} (\beta_{1} + \beta_{2}) = \lambda_{1} \beta_{1} + \lambda_{2} \beta_{2}.
\end{equation}
Notice that $ \beta_{3} = \beta_{1} + \beta_{2} $.
Collecting terms, we then have:
\begin{equation}
(s_{1} + s_{3}) \beta_{1} + (s_{2} + s_{3}) \beta_{2} = \lambda_{1} \beta_{1} + \lambda_{2} \beta_{2}.
\end{equation}
Thus we have the following linear system:
\begin{align}
s_{1} + s_{3} &= \lambda_{1}, \\
s_{2} + s_{3} &= \lambda_{2}.
\end{align}
The solution set $ S $ is:
\begin{align}
S &= \left\lbrace (\lambda_{1} - s_{3}, \lambda_{2} - s_{3}, s_{3}) \  : \  
\lambda_{1} - s_{3} \geq 0, \  \lambda_{2} - s_{3} \geq 0, \  s_{3} \geq 0 \right\rbrace  \\
&= \left\lbrace (\lambda_{1} - s_{3}, \lambda_{2} - s_{3}, s_{3}) \  : \  
s_{3} \leq \lambda_{1}, \  s_{3} \leq \lambda_{2}, \  s_{3} \geq 0 \right\rbrace  \\
&= \left\lbrace (\lambda_{1} - s_{3}, \lambda_{2} - s_{3}, s_{3}) \  : \  
0 \leq s_{3} \leq \min(\lambda_{1}, \lambda_{2}) \right\rbrace .
\end{align}
Notice that $ S = \emptyset $ if $ \lambda_{1} < 0 $ or $ \lambda_{2} < 0 $.
Therefore,
\begin{equation}\label{eq:Hbeta.lambda}
H_{\bar{\beta}}(\lambda_{1} \beta_{1} + \lambda_{2} \beta_{2}) = \vol(S) =  \cdot \max(\min(\lambda_{1}, \lambda_{2}), 0).
\end{equation}

Now we compute $ H_{\bar{\beta}}(\xi) $ for $ \xi = \ell \cdot \Omega_{1} + m \cdot (\Omega_{2} - \Omega_{1}) $.
First, we rewrite $ \xi $ in the form of $ \lambda_{1} \beta_{1} + \lambda_{2} \beta_{2} $:
\begin{align}
\xi &= \ell \cdot \Omega_{1} + m \cdot (\Omega_{2} - \Omega_{1}) \\
&= (\ell - m) \cdot \Omega_{1} + m \cdot \Omega_{2} \\
&= (\ell - m) \cdot \frac{2 \beta_{1} + \beta_{2}}{3} + m \cdot \frac{\beta_{1} + 2 \beta_{2}}{3} \\
&= \frac{2\ell - m}{3} \cdot \beta_{1} + \frac{\ell + m}{3} \cdot \beta_{2}.
\end{align}

Therefore
\begin{equation}\label{eq:trp.Hbeta}
H_{\bar{\beta}}(\xi) =  C \cdot \max\left( \min\left( \frac{2\ell - m}{3}, \frac{\ell + m}{3}\right) , 0\right) 
\end{equation}
where $C$ is a real constant.

Now consider the triple reduced product $ M_{\red}(\vecxi) $ of $ G = \SU(3) $ with the quantities $ \vecxi = (\xi_{1}, \xi_{2}, \xi_{3}) $
where
\begin{equation}
\xi_{i} = (\ell_{i} - m_{i}) \cdot \Omega_{1} + m_{i} \cdot \Omega_{2} = \ell_{i} \cdot \Omega_{1} + 
m_{i} \cdot (\Omega_{2} - \Omega_{1}), \  \ell_{i} > m_{i} > 0.
\end{equation}
Notice that each $ \xi_{i} $ lies in the open positive Weyl chamber $ \opwc $.
Our goal is to express the symplectic volume of $ M_{\red}(\vecxi) $ in terms of $ \ell_{1}, \ell_{2},
\ell_{3}, m_{1}, m_{2}, m_{3} $.

First we need to write each Weyl group element as a $ 2 \times 2 $ matrix with respect to the basis
$ \left\lbrace \Omega_{1}, \Omega_{2} - \Omega_{1} \right\rbrace  $.
Recall that we have enumerated the Weyl group as $ \left\lbrace \ws_{0}, \ws_{1}, \ws_{2}, \ws_{3}, \ws_{4},
\ws_{5} \right\rbrace  $ and associated the $ \ws_{j} $ with elements in $ \mathfrak{S}_{3} $ bijectively
through Equations (\ref{eq:wsbegin}--\ref{eq:wsend}). Each $ \ws_{j} $ can be regarded as a linear 
transformation from $ \liet $ to itself. With respect to the basis $ \left\lbrace \Omega_{1}, \Omega_{2} - \Omega_{1} \right\rbrace  $, they can be written as the following matrices:
\begin{equation} \label{weylmatrix}
\ws_{0} =
\begin{pmatrix}
1 & 0 \\
0 & 1
\end{pmatrix},
\end{equation}
\begin{equation}
\ws_{1} = 
\begin{pmatrix}
0 & 1 \\
1 & 0
\end{pmatrix},
\end{equation}
\begin{equation}
\ws_{2} =
\begin{pmatrix}
0 & -1 \\
1 & -1
\end{pmatrix},
\end{equation}
\begin{equation}
\ws_{3} =
\begin{pmatrix}
-1 & 0 \\
-1 & 1
\end{pmatrix},
\end{equation}
\begin{equation}
\ws_{4} =
\begin{pmatrix}
-1 & 1 \\
-1 & 0
\end{pmatrix},
\end{equation}
\begin{equation}
\ws_{5} =
\begin{pmatrix}
1 & -1 \\
0 & -1
\end{pmatrix}.
\end{equation}

Now, combining Equations (\ref{eq:trp.res3}) and (\ref{eq:trp.Hbeta}), we obtain our volume formula.
\end{proof}

\begin{remark}
Notice that the above formula is a piecewise linear function. While this function
is continuous, there are places, called ``walls'', where this function is not differentiable.
These walls are introduced by the $ \max $ and $ \min $ operators. If we cross a wall,
we will see that the gradient vector ``jumps''. For example,
a wall occurs when we cross the places where $ P_{ijk}^{(1)}(\vecxi) = P_{ijk}^{(2)}(\vecxi) > 0 $
for some $ i, j, k $.
\end{remark}

\section{A Result of Suzuki and Takakura}\label{sec:ST}

It will be interesting to compare our result with a result
on symplectic volume of $ N $-fold reduced products of $ G = \SU(3) $
by Suzuki and Takakura \cite{ST08} in 2008.
In this section, we will describe their result in the case when $ N = 3 $.
The settings in their paper \cite{ST08} are almost the same as ours
except that their choice of inner product on $ \lieg $ is
$ (\cdot, \cdot) / (4 \pi^{2}) $ where $ (\cdot, \cdot) $ is our choice of
inner product. Notice that this difference will not affect the symplectic
volume of triple reduced products since the symplectic volume of a coadjoint orbit 
does not depend on a choice of inner products on $ \lieg $.

Their initial input is more restrictive than ours in the following sense.
Let $ \vecxi = (\xi_{1}, \xi_{2}, \xi_{3}) $ be the quantities such that
\begin{equation}
\xi_{i} = (\ell_{i} - m_{i}) \cdot \Omega_{1} + m_{i} \cdot \Omega_{2},
\end{equation}
where $ \ell_{i} > m_{i} > 0 $ and all $ \ell_{i} $ and $ m_{i} $ are integers that
are divisible by $ 3 $ and also
\begin{equation}
(\vecw \odot \vecxi, \Omega_{1}) \neq 0
\end{equation}
for all $ \vecw \in W^{3} $.
Note that the notation $\vecw \odot \vecxi$ was introduced in (\ref{circledef}).

\begin{remark}
Notice that in our formula we do not require the $ \ell_{i} $ and the $ m_{i} $
to be integral or even rational. Thus our result is an extension of the result of
Suzuki and Takakura \cite{ST08}.
\end{remark}

Now, let $ L = \ell_{1} + \ell_{2} + \ell_{3} $ and $ M = m_{1} + m_{2} + m_{3} $.
If $ I $ is a subset of $ \left\lbrace 1, 2, 3 \right\rbrace  $,
\begin{align}
\ell_{I} &:= \sum_{i \in I} \ell_{i}, \\
m_{I} &:= \sum_{i \in I} m_{i}.
\end{align}
If $ I $ and $ J $ are two disjoint subsets of $ \left\lbrace 1, 2, 3 \right\rbrace  $,
\begin{align}
\ell_{I, J} &:= \ell_{I} + \ell_{J} = \sum_{i \in I \sqcup J} \ell_{i}, \\
m_{I, J} &:= m_{I} + m_{J} = \sum_{i \in I \sqcup J} m_{i}.
\end{align}

If $ (I_{1}, \dots, I_{6}) $ is a $ 6 $-tuple of subsets of $ \left\lbrace 1, 2, 3 \right\rbrace  $,
$ (I_{1}, \dots, I_{6}) $ is called a $ 6 $-partition of $ \left\lbrace 1, 2, 3 \right\rbrace  $
if and only if:
\begin{equation}
I_{1} \cup \dots \cup I_{6} = \left\lbrace 1, 2, 3 \right\rbrace \  \text{and} \  
I_{j} \cap I_{k} = \emptyset \  \text{whenever} \  j \neq k.
\end{equation}

For $ N = 3 $, the result of Suzuki and Takakura is the following (recall that the group $ G $
here is $ \SU(3) $):

\begin{theorem}[Theorem 4.5 in \cite{ST08}, in the case $ N = 3 $]\label{thm:ST08.N=3}
Let $ \mathcal{I}_{\vecxi} $ denote the set of those $ 6 $-partitions $ (I_{1}, \dots, I_{6}) $
of $ \left\lbrace 1, 2, 3 \right\rbrace  $ such that
\begin{align}
\ell_{I_{1}, I_{2}} + m_{I_{4}, I_{5}} &< \frac{L + M}{3}, \  \text{and} \\
\ell_{I_{3}, I_{4}} + m_{I_{6}, I_{1}} &< \frac{L + M}{3}.
\end{align}
Let $ \mathcal{J}_{\vecxi} $ denote the set of those $ 6 $-partitions $ (I_{1}, \dots, I_{6}) $
of $ \left\lbrace 1, 2, 3 \right\rbrace  $ such that
\begin{align}
\ell_{I_{3}, I_{4}} + m_{I_{6}, I_{1}} &> \frac{L + M}{3}, \  \text{and} \\
\ell_{I_{5}, I_{6}} + m_{I_{2}, I_{3}} &> \frac{L + M}{3}.
\end{align}
Let $ A_{\vecxi}: \mathcal{I}_{\vecxi} \to \RR $ be defined by:
\begin{equation}
A_{\vecxi}(I_{1}, \dots, I_{6}) := \frac{-(-1)^{\left| I_{1} \right| + \left| I_{3} \right| + \left| I_{5} \right| }}{6} \left( \frac{L + M}{3} - \ell_{I_{1}, I_{2}} - m_{I_{4}, I_{5}} \right) .
\end{equation}
Let $ B_{\vecxi}: \mathcal{J}_{\vecxi} \to \RR $ be defined by:
\begin{equation}
B_{\vecxi}(I_{1}, \dots, I_{6}) := \frac{-(-1)^{\left| I_{1} \right| + \left| I_{3} \right| + \left| I_{5} \right| }}{6} \left( \ell_{I_{5}, I_{6}} + m_{I_{2}, I_{3}} - \frac{L + M}{3} \right) .
\end{equation}
Then, the symplectic volume of $ M_{\red}(\vecxi) $ is given by:
\begin{equation}\label{eq:trp.vol.ST08}
\mathcal{V}(\vecxi) = \sum_{(I_{1}, \dots, I_{6}) \in \mathcal{I}_{\vecxi}} A_{\vecxi}(I_{1}, \dots, I_{6})
+ \sum_{(I_{1}, \dots, I_{6}) \in \mathcal{J}_{\vecxi}} B_{\vecxi}(I_{1}, \dots, I_{6}).
\end{equation}
\end{theorem}

We shall briefly explain this result.

First, $ 6 $-partitions $ (I_{1}, \dots, I_{6}) $ of $ \left\lbrace 1, 2, 3 \right\rbrace  $
correspond bijectively to fixed points of $ T $ acting on $ M $, i.e., to $ M^{T} $, in the following way.
Suzuki and Takakura have used a different enumeration
$ \left\lbrace \sigma_{1}, \sigma_{2}, \sigma_{3}, \sigma_{4}, \sigma_{5}, \sigma_{6} \right\rbrace  $ of the Weyl group $ W $ of $ \SU(3) $ in their paper \cite{ST08}.
Relating their enumeration with ours, we have:
\begin{align}
\sigma_{1} &= \Id = \ws_{0}, \\
\sigma_{2} &= (2 \  3) = \ws_{5}, \\
\sigma_{3} &= (1 \  2 \  3) = \ws_{2}, \\
\sigma_{4} &= (1 \  2) = \ws_{1}, \\
\sigma_{5} &= (1 \  3 \  2) = \ws_{4}, \\
\sigma_{6} &= (1 \  3) = \ws_{3}.
\end{align}

For each $ \vecw = (w_{1}, w_{2}, w_{3}) \in W^{3} $,
let $ I_{j} $ be defined as
\begin{equation}
I_{j} := \left\lbrace i \in \left\lbrace 1, 2, 3 \right\rbrace  \  : \  w_{i} = \sigma_{j} \right\rbrace 
\end{equation}
for $ j = 1, \dots, 6 $.
Then $ (I_{1}, \dots, I_{6}) $ is a $ 6 $-partition of $ \left\lbrace 1, 2, 3 \right\rbrace  $.
For example, given $ (\sigma_{2}, \sigma_{5}, \sigma_{2}) \in W^{3} $, the corresponding $ 6 $-partition
of $ \left\lbrace 1, 2, 3 \right\rbrace  $ is
\begin{equation}
\left( \emptyset, \left\lbrace 1, 3 \right\rbrace , \emptyset, \emptyset, \left\lbrace 2 \right\rbrace , \emptyset \right) .
\end{equation}
In other words, $ I_{j} $ tells us in which coordinates (in this case, the first, the second, or the third
coordinate of $ \vecw $) $ \sigma_{j} $ appears in $ \vecw $.

Suzuki and Takakura have observed that $ \vecw \odot \vecxi $ is the matrix
\begin{equation}
2 \pi \I \cdot \diag(\ell_{I_{1}, I_{2}} + m_{I_{4}, I_{5}} - \frac{L + M}{3},
\ell_{I_{3}, I_{4}} + m_{I_{6}, I_{1}} - \frac{L + M}{3},
\ell_{I_{5}, I_{6}} + m_{I_{2}, I_{3}} - \frac{L + M}{3}).
\end{equation}
We would like to determine under what conditions the above matrix is in the cone spanned by $ \bar{\beta} =
(\beta_{1}, \beta_{2}, \beta_{3}) $.
The vector $ \vecw \odot \vecxi $ is in the cone spanned by $ \bar{\beta} $ if and only if
\begin{align}
(\vecw \odot \vecxi, \Omega_{1}) &> 0, \  \text{and} \\
(\vecw \odot \vecxi, \Omega_{2}) &> 0.
\end{align}
Recall that
\begin{align}
\Omega_{1} &= \frac{2 \pi \I}{3} \cdot \diag(2, -1, -1), \\
\Omega_{2} &= \frac{2 \pi \I}{3} \cdot \diag(1, 1, -2).
\end{align}
Thus, $ \vecw \odot \vecxi $ is in the cone spanned by $ \bar{\beta} $ if and only if
\begin{align}
2 \ell_{I_{1}, I_{2}} + 2 m_{I_{4}, I_{5}} - \ell_{I_{3}, I_{4}} - m_{I_{6}, I_{1}}
- \ell_{I_{5}, I_{6}} - m_{I_{2}, I_{3}} &> 0, \  \text{and} \\
\ell_{I_{1}, I_{2}} + m_{I_{4}, I_{5}} + \ell_{I_{3}, I_{4}} + m_{I_{6}, I_{1}}
- 2 \ell_{I_{5}, I_{6}} - 2 m_{I_{2}, I_{3}} &> 0.
\end{align}
Notice that
\begin{align}
(\vecw \odot \vecxi, \Omega_{1}) &= \frac{4 \pi^{2}}{3}
\left( 2 \ell_{I_{1}, I_{2}} + 2 m_{I_{4}, I_{5}} - \ell_{I_{3}, I_{4}} - m_{I_{6}, I_{1}}
- \ell_{I_{5}, I_{6}} - m_{I_{2}, I_{3}} \right)  \\
&= 4 \pi^{2} \left( \ell_{I_{1}, I_{2}} + m_{I_{4}, I_{5}} - \frac{L + M}{3} \right) , \  \text{and} \\
(\vecw \odot \vecxi, \Omega_{2}) &= \frac{4 \pi^{2}}{3}
\left( \ell_{I_{1}, I_{2}} + m_{I_{4}, I_{5}} + \ell_{I_{3}, I_{4}} + m_{I_{6}, I_{1}}
- 2 \ell_{I_{5}, I_{6}} - 2 m_{I_{2}, I_{3}} \right)  \\
&= 4 \pi^{2} \left( \frac{L + M}{3} - \ell_{I_{5}, I_{6}} - m_{I_{2}, I_{3}} \right) .
\end{align}

It will be interesting to compare the above condition with the conditions for $ \mathcal{I}_{\vecxi} $
and $ \mathcal{J}_{\vecxi} $ used in Theorem \ref{thm:ST08.N=3}.
If $ \vecw \odot \vecxi $ satisfies the condition for $ \mathcal{I}_{\vecxi} $, then
\begin{align}
&2 \ell_{I_{1}, I_{2}} + 2 m_{I_{4}, I_{5}} - \ell_{I_{3}, I_{4}} - m_{I_{6}, I_{1}}
- \ell_{I_{5}, I_{6}} - m_{I_{2}, I_{3}} \\
= \  &3 \ell_{I_{1}, I_{2}} + 3 m_{I_{4}, I_{5}} - L - M \\
= \  &3 (\ell_{I_{1}, I_{2}} + m_{I_{4}, I_{5}} - \frac{L + M}{3}) < 0.
\end{align}
Thus, $ \vecw \odot \vecxi $ is not in the cone spanned by $ \bar{\beta} $.

If we look closely at the first inequality of the condition for $ \mathcal{I}_{\vecxi} $, namely,
\begin{equation}
\ell_{I_{1}, I_{2}} + m_{I_{4}, I_{5}} < \frac{L + M}{3},
\end{equation}
this means exactly
\begin{equation}
(\vecw \odot \vecxi, \Omega_{1}) < 0.
\end{equation}
Also, the second inequality of the condition for $ \mathcal{I}_{\vecxi} $, namely,
\begin{equation}
\ell_{I_{3}, I_{4}} + m_{I_{6}, I_{1}} < \frac{L + M}{3},
\end{equation}
means exactly
\begin{equation}
(\vecw \odot \vecxi, \Omega_{2} - \Omega_{1}) < 0.
\end{equation}
Thus, we have translated the condition for $ \mathcal{I}_{\vecxi} $ into the following two inequalities:
\begin{align}
(\vecw \odot \vecxi, \Omega_{1}) &< 0, \  \text{and} \\
(\vecw \odot \vecxi, \Omega_{2} - \Omega_{1}) &< 0.
\end{align}

Similarly, we can translate the condition for $ \mathcal{J}_{\vecxi} $ into the following two inequalities:
\begin{align}
(\vecw \odot \vecxi, \Omega_{2} - \Omega_{1}) &> 0, \  \text{and} \\
(\vecw \odot \vecxi, \Omega_{2}) &< 0.
\end{align}

Notice that if $ \vecw \odot \vecxi $ satisfies either the condition for $ \mathcal{I}_{\vecxi} $ or
the condition for $ \mathcal{J}_{\vecxi} $, we always have that $ \vecw \odot \vecxi $ is in the cone
spanned by $ -\bar{\beta} = (-\beta_{1}, -\beta_{2}, -\beta_{1} - \beta_{2}) $. In other words,
only those $ \vecw \odot \vecxi $ contained in the cone spanned by $ -\bar{\beta} $ will contribute to
the sum in the volume formula of Suzuki and Takakura. Also notice that the sign
\begin{equation}
-(-1)^{\left| I_{1} \right| + \left| I_{3} \right| + \left| I_{5} \right| }
\end{equation}
is exactly the signature of $ \vecw $.

The above observations lead us to conclude that our volume formula (Theorem \ref{thm:trp.vol}) and the volume
formula of Suzuki and Takakura (Theorem \ref{thm:ST08.N=3}) are very closely related.

We have the following result.

\begin{theorem}\label{thm:comparison.trp}
Under the assumptions (A1), (A2) and (A3) stated in \S \ref{sss:assumptions}, our volume formula in Theorem \ref{thm:trp.vol} 
agrees completely with  the volume formula of
Suzuki and Takakura (Theorem \ref{thm:ST08.N=3}) for triple reduced products of $ \SU(3) $, provided that
$ K = -K^{\prime} = - 1/6 $. Thus our volume formula extends that of 
\cite{ST08}.
\end{theorem}

\begin{proof}
The first key observation supporting that our formula should indeed agree with theirs is that
we derived our formula by using the residue formula (Theorem \ref{thm:residue}) with a choice of
cone, namely the cone spanned by $ \bar{\beta} $ and as a result, only those $ \vecw \odot \vecxi $
in the cone spanned by $ \bar{\beta} $ will contribute to the sum in our volume formula. However,
the total sum in the residue formula does not depend on the choice of cone. Therefore, we could equally
well choose the cone spanned by $ -\bar{\beta} $ to carry out the computations of the individual terms
in the sum. Let us  carry this out.

More precisely, we start from the choice of $ \Lambda_{-} = -\opwc $. In this way, all of 
$ -\beta_{1}, -\beta_{2}, -\beta_{3} = -\beta_{1} - \beta_{2} $ lie in the dual cone $ \Lambda_{-}^{*} $.

With this new choice of cone, we carry out the computation, starting from Equation (\ref{eq:res.totalsum.N=3}):
\begin{align}
&\res\left( \sum_{\vecw \in W^{3}} \sgn(\vecw) \frac{e^{\I (\vecw \odot \vecxi, \psi)} [d\psi]}{\varpi(\psi)} \right) \\
= &\res^{\Lambda_{-}}\left( \sum_{\vecw \in W^{3}} \sgn(\vecw) \frac{e^{\I (\vecw \odot \vecxi, \psi)} [d\psi]}{\prod_{j = 1}^{3} \beta_{j}(\psi)} \right) \\
= &\sum_{\vecw \in W^{3}} \sgn(\vecw) \res^{\Lambda_{-}}\left( \frac{e^{\I (\vecw \odot \vecxi, \psi)} [d\psi]}{\prod_{j = 1}^{3} \beta_{j}(\psi)} \right) \\
= &-\sum_{\vecw \in W^{3}} \sgn(\vecw) \res^{\Lambda_{-}}\left( \frac{e^{\I (\vecw \odot \vecxi, \psi)} [d\psi]}{\prod_{j = 1}^{3} (-\beta_{j})(\psi)} \right) \\
= &-\sum_{\vecw \in W^{3}} \sgn(\vecw) \frac{\I^{3}}{(2 \pi \I)^{2}} H_{(-\bar{\beta})}(\vecw \odot \vecxi).
\end{align}

Now our volume formula corresponding to the cone $ \Lambda_{-} $ can be written as
\begin{equation}
\vol^{\SSS}(M_{\red}(\vecxi)) = K^{\prime} \sum_{\vecw \in W^{3}} \sgn(\vecw) H_{(-\bar{\beta})}(\vecw \odot \vecxi),
\end{equation}
where $ K^{\prime} $ is a  constant. Notice that the constant $ K $ in our volume formula corresponding to the cone $ \Lambda $, i.e., Equation (\ref{eq:trp.vol}), is simply  
$ K = -K^{\prime}$.

The function $ H_{(-\bar{\beta})} $ is supported in the cone spanned by $ -\bar{\beta} $.
Therefore, only those $ \vecw \odot \vecxi $ inside this cone will contribute to the sum. This
gives us the common ground to compare our volume formula (using the cone $ \Lambda_{-} $) and
the volume formula of Suzuki and Takakura.

Let us  compute $ H_{(-\bar{\beta})}(\lambda_{1} \cdot (-\beta_{1}) + \lambda_{2} \cdot (-\beta_{2})) $ where $ \lambda_{1}, \lambda_{2} \in \RR $.

We have:
\begin{align}
&H_{(-\bar{\beta})}(\lambda_{1} \cdot (-\beta_{1}) + \lambda_{2} \cdot (-\beta_{2})) \\
= &\vol\left\lbrace (s_{1}, s_{2}, s_{3}) \in \RR_{+}^{3} \  : \  
\sum_{j = 1}^{3} s_{j} \cdot (-\beta_{j}) = \lambda_{1} \cdot (-\beta_{1}) + \lambda_{2} \cdot (-\beta_{2}) \right\rbrace .
\end{align}
Thus, we need to solve the following equation:
\begin{equation}
s_{1} \cdot (-\beta_{1}) + s_{2} \cdot (-\beta_{2}) + s_{3} \cdot (-\beta_{1} - \beta_{2})
= \lambda_{1} \cdot (-\beta_{1}) + \lambda_{2} \cdot (-\beta_{2}).
\end{equation}
Therefore, we need to solve the following linear system:
\begin{align}
s_{1} + s_{3} &= \lambda_{1}, \\
s_{2} + s_{3} &= \lambda_{2}.
\end{align}
The solution set $ S_{-} $ is:
\begin{align}
S_{-} &= \left\lbrace (\lambda_{1} - s_{3}, \lambda_{2} - s_{3}, s_{3}) \  : \  
\lambda_{1} - s_{3} \geq 0, \  \lambda_{2} - s_{3} \geq 0, \  s_{3} \geq 0 \right\rbrace \\
&= \left\lbrace (\lambda_{1} - s_{3}, \lambda_{2} - s_{3}, s_{3}) \  : \  
s_{3} \leq \lambda_{1}, \  s_{3} \leq \lambda_{2}, s_{3} \geq 0 \right\rbrace \\
&= \left\lbrace (\lambda_{1} - s_{3}, \lambda_{2} - s_{3}, s_{3}) \  : \  
0 \leq s_{3} \leq \min(\lambda_{1}, \lambda_{2}) \right\rbrace .
\end{align}
Therefore,
\begin{equation}
H_{-\bar{\beta}}(\lambda_{1} \cdot (-\beta_{1}) + \lambda_{2} \cdot (-\beta_{2})) = \vol(S_{-})
= C \cdot \max(\min(\lambda_{1}, \lambda_{2}), 0),
\end{equation}
where $ C $ is the same constant as in Equation (\ref{eq:Hbeta.lambda}).
(In fact $C = 1.$ )

Now, let's look at those $ \vecw \odot \vecxi $ inside the cone spanned by $ -\bar{\beta} $.
Without loss of generality we can assume that $ \vecxi $ is generic so that for all $ \vecw \in W^{3} $,
\begin{align}
(\vecw \odot \vecxi, \Omega_{1}) &\neq 0, \\
(\vecw \odot \vecxi, \Omega_{2}) &\neq 0, \\
(\vecw \odot \vecxi, \Omega_{2} - \Omega_{1}) &\neq 0.
\end{align}
Therefore, the collection of those $ \vecw \odot \vecxi $ inside the cone spanned by $ -\bar{\beta} $ is
the disjoint union of the two sets $ \mathcal{A}_{\vecxi} $ and $ \mathcal{B}_{\vecxi} $, where
$ \mathcal{A}_{\vecxi} $ denotes the set of those $ \vecw \odot \vecxi $ such that
\begin{align}
(\vecw \odot \vecxi, \Omega_{1}) &< 0, \  \text{and} \\
(\vecw \odot \vecxi, \Omega_{2} - \Omega_{1}) &< 0,
\end{align}
and $ \mathcal{B}_{\vecxi} $ denotes the set of those $ \vecw \odot \vecxi $ such that
\begin{align}
(\vecw \odot \vecxi, \Omega_{2} - \Omega_{1}) &> 0, \  \text{and} \\
(\vecw \odot \vecxi, \Omega_{2}) &< 0.
\end{align}

Notice that the above grouping is in complete agreement with the grouping by $ \mathcal{I}_{\vecxi} $ and
$ \mathcal{J}_{\vecxi} $ in the formula of Suzuki and Takakura so that we can make a term-by-term
comparison between our formula and theirs.

For each $ \vecw \odot \vecxi \in \mathcal{A}_{\vecxi} $, it is easy to see that
$ \vecw \odot \vecxi $ can be written as $ \lambda_{1} \cdot (-\beta_{1}) + \lambda_{2} \cdot (-\beta_{2})$ with some $ \lambda_{1}, \lambda_{2} $ satisfying $ 0 < \lambda_{1} < \lambda_{2} $.
Therefore, the contribution of this $ \vecw \odot \vecxi $ to our volume formula is
\begin{equation}
K^{\prime} \cdot \sgn(\vecw) \cdot H_{(-\bar{\beta})}(\vecw \odot \vecxi) =
K^{\prime} \cdot \sgn(\vecw) \cdot \lambda_{1}.
\end{equation}

Notice that this $ \vecw $ will correspond to a $ 6 $-partition $ (I_{1}, \dots, I_{6}) $ and we have
\begin{equation}
\sgn(\vecw) = -(-1)^{\left| I_{1} \right| + \left| I_{3} \right| + \left| I_{5} \right| }.
\end{equation}
Now we only need to figure out how we can express this $ \lambda_{1} $ in terms of the $ \ell_{i} $ and 
$ m_{i} $.

First, given $ a, b \in \RR $ and 
\begin{equation}
2 \pi \I \cdot \diag(a, b, -a - b) = \lambda_{1} \cdot (-\beta_{1}) + \lambda_{2} \cdot (-\beta_{2}),
\end{equation}
we want to express $ \lambda_{1} $ and $ \lambda_{2} $ in terms of $ a, b $. This is equivalent to
solving the following linear system:
\begin{align}
-\lambda_{1} &= a, \\
\lambda_{1} - \lambda_{2} &= b.
\end{align}
Thus, we have
\begin{align}
\lambda_{1} &= -a, \\
\lambda_{2} &= -a - b.
\end{align}

Recall that $ \vecw \odot \vecxi $ is the matrix
\begin{equation}
2 \pi \I \cdot \diag(\ell_{I_{1}, I_{2}} + m_{I_{4}, I_{5}} - \frac{L + M}{3},
\ell_{I_{3}, I_{4}} + m_{I_{6}, I_{1}} - \frac{L + M}{3},
\ell_{I_{5}, I_{6}} + m_{I_{2}, I_{3}} - \frac{L + M}{3}).
\end{equation}

Hence, if $ \vecw \odot \vecxi \in \mathcal{A}_{\vecxi} $,
the contribution of this $ \vecw \odot \vecxi $ to our volume formula is
\begin{equation}
K^{\prime} \cdot \sgn(\vecw) \cdot \lambda_{1} =
K^{\prime} \cdot \left( -(-1)^{\left| I_{1} \right| + \left| I_{3} \right| + \left| I_{5} \right| } \right)  \cdot \left( \frac{L + M}{3} - \ell_{I_{1}, I_{2}} - m_{I_{4}, I_{5}} \right) ,
\end{equation}
which precisely matches the term
\begin{equation}
A_{\vecxi}(I_{1}, \dots, I_{6}) = \frac{-(-1)^{\left| I_{1} \right| + \left| I_{3} \right| + \left| I_{5} \right| }}{6} \left( \frac{L + M}{3} - \ell_{I_{1}, I_{2}} - m_{I_{4}, I_{5}} \right) 
\end{equation}
for the contribution of this $ \vecw \odot \vecxi $ to the volume formula of Suzuki and Takakura, provided that $ K^{\prime} = 1/6 $.

For each $ \vecw \odot \vecxi \in \mathcal{B}_{\vecxi} $, it is easy to see that
$ \vecw \odot \vecxi $ can be written as $ \lambda_{1} \cdot (-\beta_{1}) + \lambda_{2} \cdot (-\beta_{2})$ with some $ \lambda_{1}, \lambda_{2} $ satisfying $ 0 < \lambda_{2} < \lambda_{1} $.
Therefore, the contribution of this $ \vecw \odot \vecxi $ to our volume formula is
\begin{align}
&K^{\prime} \cdot \sgn(\vecw) \cdot H_{(-\bar{\beta})}(\vecw \odot \vecxi) \\
= &K^{\prime} \cdot \sgn(\vecw) \cdot \lambda_{2} \\
= &K^{\prime} \cdot \left( -(-1)^{\left| I_{1} \right| + \left| I_{3} \right| + \left| I_{5} \right| } \right) \cdot \left( \ell_{I_{5}, I_{6}} + m_{I_{2}, I_{3}} - \frac{L + M}{3} \right) ,
\end{align}
which precisely matches the term
\begin{equation}
B_{\vecxi}(I_{1}, \dots, I_{6}) = \frac{-(-1)^{\left| I_{1} \right| + \left| I_{3} \right| + \left| I_{5} \right| }}{6} \left( \ell_{I_{5}, I_{6}} + m_{I_{2}, I_{3}} - \frac{L + M}{3} \right) 
\end{equation}
for the contribution of this $ \vecw \odot \vecxi $ to the volume formula of Suzuki and Takakura, provided that $ K^{\prime} = 1/6 $.

Observing  that the  sum in the residue formula does not depend on the choice of cone,
we have proved the theorem.
\end{proof}

\section{Generalizations of Volume Formula}\label{ch:generalizations}

In this section, we generalize some of our earlier results.

\subsection{Volume Formula for $ N $-fold Reduced Products of $ \SU(3) $}

In this section, our group $ G $ is still $ \SU(3) $. As before, let $ T $ be the standard
maximal torus in $ G $. Let $ W $ denote the Weyl group.
Thus, $ W = \mathfrak{S}_{3} $.

We assume the following.
\begin{itemize}
\item  $ N \geq 3 $ is a positive integer. Notice that previously
our $ N $ was equal to $ 3 $.

\item Suppose  $ \vecxi = (\xi_{1}, \dots, \xi_{N}) $ is a collection of $N$ elements in 
$\liet_+$ which 
 satisfy the conditions specified in Section \ref{ch:trp}.
\item Let $ M = \orb_{\xi_{1}} \times \dots \times \orb_{\xi_{N}} $. This is a compact
symplectic manifold with $ G $ acting diagonally on it in a Hamiltonian fashion.
\item Let $ \mu_{G} $ and $ \mu_{T} $ be the moment maps of the $ G $-action and the $ T $-action
respectively.
\item Let $ M^{T}, M_{0}, M_{\red} $ be defined similarly as in \S \ref{ch:trp}.
\end{itemize}
The following proposition is an easy generalization of Proposition \ref{prop:TfixedPoints.N=3}.

\begin{proposition}\label{prop:TfixedPoints.N}
Let $ G = \SU(3) $ and $ M = \orb_{\xi_{1}} \times \dots \times \orb_{\xi_{N}} $ be
the Cartesian product of $ N \geq 3 $ adjoint orbits of $ G $, where the $ \xi_{i} $
satisfy the assumptions (A1) and (A2). Then,
$ M^{T} $ is the discrete set
\begin{equation}
\left\lbrace (w_{1} \cdot \xi_{1}, \dots, w_{N} \cdot \xi_{N}) \  : \  
w_{i} \in W \right\rbrace .
\end{equation}
Thus, $ \left| M^{T} \right| = \left| W \right| ^{N} $.
\end{proposition}

As a result, $ M^{T} $ is parametrized by $ \vecw \in W^{N} $. Moreover, $ 6 $-partitions of
$ \left\lbrace 1, \dots, N \right\rbrace  $ can be defined similarly and each $ 6 $-partition $ (I_{1}, \dots, I_{6}) $ of
$ \left\lbrace 1, \dots, N \right\rbrace  $ corresponds to a unique $ \vecw \in W^{N} $ in the same way
as before:
\begin{equation}
I_{j} = \left\lbrace i \in \left\lbrace 1, \dots, N \right\rbrace  \  : \  w_{i} = \sigma_{j} \right\rbrace ,
\end{equation}
for $ j = 1, \dots, 6 $. Thus, $ M^{T} $ can also be parametrized by all the $ 6 $-partitions of
$ \left\lbrace 1, \dots, N \right\rbrace  $.

In addition, $ \vecw \cdot \vecxi $ and $ \vecw \odot \vecxi $ can similarly be  defined. Note that
the notation $\vecw \odot \vecxi$ was introduced earlier in (\ref{circledef}).

We have the following result, whose proof is similar to the  proof
of our earlier result  \ref{thm:trp.vol}): 

\begin{theorem}\label{thm:Nfold.vol}
Let $ G = \SU(3) $ and $ M = \orb_{\xi_{1}} \times \dots \times \orb_{\xi_{N}} $ be
the Cartesian product of $ N \geq 3 $ adjoint orbits of $ G $, where the $ \xi_{i} $'s
satisfy the assumptions (A1), (A2) and (A3). 
Here $s$ is the real dimension of $G$ and $l$ is the real dimension of $T$
In this case $G = SU(3)$, so $s = 8 $ and $l = 2$. Therefore
$d = 6N-16.$
Then,
the symplectic volume of $ M_{\red}(\vecxi) $ is:
\begin{equation}\label{eq:Nfold.vol}
\vol^{\SSS}(M_{\red}(\vecxi)) = \frac{1}{\I^{d/2}} n_{0} C_{G} 
\res\left( \sum_{\vecw \in W^{N}} \sgn(\vecw) \frac{e^{\I (\vecw \odot \vecxi, \psi)} [d\psi]}{\varpi^{N - 2}(\psi)} \right) ,
\end{equation}
where $ n_{0} $ and $ C_{G} $ are as the same as in Theorem \ref{thm:residue} and
\begin{equation}
d = N (s - l) - 2s = (N - 2) s - N l,
\end{equation}
where $ s $ is the real dimension of $ G $ and $ l $ is the real dimension of $ T $.

\end{theorem}

\begin{proof}
When computing the symplectic volume of $ M_{\red}(\vecxi) $, the only essential difference between
the $ N = 3 $ case and the general $ N \geq 3 $ case lies in the equivariant Euler class $ \ee_{F} $
of the normal bundle of the fixed points of $ T $.
In the general $ N \geq 3 $ case, for each fixed point $ F = \vecw \cdot \vecxi \in M^{T} $,
\begin{equation}
\ee_{F}(\psi) = \sgn(\vecw) \cdot \varpi^{N}(\psi).
\end{equation}
This completes the proof.
\end{proof}

As in Section \ref{sec:ST}, we should try  to compare our formula with the general formula of Suzuki and
Takakura \cite{ST08}. Therefore, by the comparison argument in Section \ref{sec:ST}, when computing the residue, we will use the cone $ \Lambda_{-} = -\opwc $ instead of $ \Lambda = \opwc $. However,
as we will see soon, the comparison for general $ N $ seems more  difficult and we do not have a complete
comparison for general $ N $ case.

We first write down here the general formula of Suzuki and Takakura
\cite{ST08} (recall that the group $ G $ here is $ \SU(3) $):

\begin{theorem}[Theorem 4.5 in \cite{ST08}]: (Compare with  (\ref{thm:ST08.N=3}), the case $N=3$.)
Let $ N \geq 3 $ be an integer. Let
\begin{equation}
\xi_{i} = (\ell_{i} - m_{i}) \cdot \Omega_{1} + m_{i} \cdot \Omega_{2},
\end{equation}
where $ \ell_{i} > m_{i} > 0 $ are all integers divisible by $ 3 $ and
\begin{equation}
(\vecw \odot \vecxi, \Omega_{1}) \neq 0
\end{equation}
for all $ \vecw \in W^{N} $. Let
\begin{equation}
L = \sum_{i = 1}^{N} \ell_{i}, \quad M = \sum_{i = 1}^{N} m_{i}.
\end{equation}
Let $ \mathcal{I}_{\vecxi} $ denote the set of $ 6 $-partitions $ (I_{1}, \dots, I_{6}) $ of
$ \left\lbrace 1, \dots, N \right\rbrace  $ such that
\begin{align}
\ell_{I_{1}, I_{2}} + m_{I_{4}, I_{5}} &< \frac{L + M}{3}, \\
\ell_{I_{3}, I_{4}} + m_{I_{6}, I_{1}} &< \frac{L + M}{3}.
\end{align}
Let $ \mathcal{J}_{\vecxi} $ denote the set of $ 6 $-partitions $ (I_{1}, \dots, I_{6}) $ of
$ \left\lbrace 1, \dots, N \right\rbrace  $ such that
\begin{align}
\ell_{I_{3}, I_{4}} + m_{I_{6}, I_{1}} &> \frac{L + M}{3}, \\
\ell_{I_{5}, I_{6}} + m_{I_{2}, I_{3}} &> \frac{L + M}{3}.
\end{align}
Let $ A_{\vecxi}: \mathcal{I}_{\vecxi} \to \RR $ be defined by
\begin{align}
A_{\vecxi}(I_{1}, \dots, I_{6}) := &\frac{-(-1)^{\left| I_{1} \right| + \left| I_{3} \right| + \left| I_{5} \right| }}{6 (3N - 8)!} \sum_{j = 0}^{N - 3} \binom{3N - 8}{j} \binom{2N - 6 - j}{N - 3} \notag \\
&\left( \frac{L + M}{3} - \ell_{I_{3}, I_{4}} - m_{I_{6}, I_{1}} \right) ^{j} \left( \frac{L + M}{3} - \ell_{I_{1}, I_{2}} - m_{I_{4}, I_{5}} \right) ^{3N - 8 - j}.
\end{align}
Let $ B_{\vecxi}: \mathcal{J}_{\vecxi} \to \RR $ be defined by
\begin{align}
B_{\vecxi}(I_{1}, \dots, I_{6}) := &\frac{-(-1)^{\left| I_{1} \right| + \left| I_{3} \right| + \left| I_{5} \right| }}{6 (3N - 8)!} \sum_{j = 0}^{N - 3} \binom{3N - 8}{j} \binom{2N - 6 - j}{N - 3} \notag \\
&\left( \ell_{I_{3}, I_{4}} + m_{I_{6}, I_{1}} - \frac{L + M}{3} \right) ^{j} \left( \ell_{I_{5}, I_{6}} + m_{I_{2}, I_{3}} - \frac{L + M}{3} \right) ^{3N - 8 - j}.
\end{align}
Then, the symplectic volume of $ M_{\red}(\vecxi) $ is given by
\begin{equation}\label{eq:ST.vol.general}
\mathcal{V}(\vecxi) = \sum_{(I_{1}, \dots, I_{6}) \in \mathcal{I}_{\vecxi}} A_{\vecxi}(I_{1}, \dots, I_{6}) + \sum_{(I_{1}, \dots, I_{6}) \in \mathcal{J}_{\vecxi}} B_{\vecxi}(I_{1}, \dots, I_{6}).
\end{equation}
\end{theorem}


We have the following result. 

\begin{theorem}
Let $ G = \SU(3) $ and $ M = \orb_{\xi_{1}} \times \dots \times \orb_{\xi_{N}} $ be
the Cartesian product of $ N \geq 3 $ adjoint orbits of $ G $, where the $ \xi_{i} $'s
satisfy the assumptions (A1), (A2) and (A3). Then,
the symplectic volume of $ M_{\red}(\vecxi) $ is:
\begin{equation}
\vol^{\SSS}(M_{\red}(\vecxi)) = C \cdot \sum_{\vecw \in W^{N}} \sgn(\vecw) H_{(-\bar{\beta})^{N-2}}(\vecw \odot \vecxi),
\end{equation}
where $ C $ is a constant.
\end{theorem}

\begin{proof}
Recall that all of $ -\beta_{1}, -\beta_{2}, -\beta_{3} = -\beta_{1} - \beta_{2} $ lie in the dual
cone $ \Lambda_{-}^{*} $.

We have:
\begin{align}
&\res\left( \sum_{\vecw \in W^{N}} \sgn(\vecw) \frac{e^{\I (\vecw \odot \vecxi, \psi)} [d\psi]}{\varpi^{N - 2}(\psi)} \right) \\
= \  &\res^{\Lambda_{-}}\left( \sum_{\vecw \in W^{N}} \sgn(\vecw) \frac{e^{\I (\vecw \odot \vecxi, \psi)} [d\psi]}{\varpi^{N - 2}(\psi)} \right) \\
= \  &\sum_{\vecw \in W^{N}} \sgn(\vecw) \res^{\Lambda_{-}}\left( \frac{e^{\I (\vecw \odot \vecxi, \psi)} [d\psi]}{\varpi^{N - 2}(\psi)} \right) \\
= \  &\sum_{\vecw \in W^{N}} \sgn(\vecw) \res^{\Lambda_{-}}\left( \frac{e^{\I (\vecw \odot \vecxi, \psi)} [d\psi]}{\left( (-1) \cdot \prod_{j = 1}^{3} (-\beta_{j})(\psi) \right) ^{N-2}} \right) \\
= \  &(-1)^{N-2} \cdot \sum_{\vecw \in W^{N}} \sgn(\vecw) \res^{\Lambda_{-}}\left( \frac{e^{\I (\vecw \odot \vecxi, \psi)} [d\psi]}{\left( \prod_{j = 1}^{3} (-\beta_{j})(\psi) \right) ^{N-2}} \right) .
\end{align}

To compute
\begin{equation}
\res^{\Lambda_{-}}\left( \frac{e^{\I (\vecw \odot \vecxi, \psi)} [d\psi]}{\left( \prod_{j = 1}^{3} (-\beta_{j})(\psi) \right) ^{N-2}} \right) ,
\end{equation}
we introduce the notation $ (-\bar{\beta})^{N-2} $ to denote the following:
\begin{equation}
(-\bar{\beta})^{N-2} := (-\beta_{1}, -\beta_{2}, -\beta_{3}, \dots, -\beta_{1}, -\beta_{2}, -\beta_{3}),
\end{equation}
where the sequence $ -\beta_{1}, -\beta_{2}, -\beta_{3} $ repeats itself for $ N - 2 $ times.

Now we have:
\begin{align}
&\res^{\Lambda_{-}}\left( \frac{e^{\I (\vecw \odot \vecxi, \psi)} [d\psi]}{\left( \prod_{j = 1}^{3} (-\beta_{j})(\psi) \right) ^{N-2}} \right) \\
= \  &\frac{\I^{3(N-2)}}{(2 \pi \I)^{2}} \cdot H_{(-\bar{\beta})^{N-2}}(\vecw \odot \vecxi).
\end{align}
This completes the proof.
\end{proof}


\subsection{Volume Formula for General $ N $-fold Reduced Products}\label{sec:generalVol}

The method of nonabelian localization and the residue formula apply not only for $ G = \SU(3) $, but
also for any compact connected Lie groups. However, to apply Theorem \ref{thm:residue} in our situation,
namely the situation where the group $ G $ acts diagonally on the product of adjoint orbits by the adjoint action, we need to
make sure that the stabilizer of any point in $ M_{0} = \mu_{G}^{-1}(0) $ is finite. Therefore, in addition
to the Lie group $ G $ being compact and connected, we assume that $ G $ is also semisimple.  

We note that in \cite{ST09}  Suzuki and Takakura also treat 
volumes of $N$-fold reduced products of compact Lie groups $G$. It 
would be interesting to verify explicitly that their results agree with ours.
Our methods  may be more amenable to generalization to 
intersection pairings, a subject we treat in \cite{JJ19}.

Let $ T $ denote a chosen maximal torus of $ G $. Let $ \lieg $ be the Lie algebra of $ G $ and $ \liet $
be the Lie algebra of $ T $. Let $ W $ denote the Weyl group $ \NN(T) / T $. Let $ (\cdot, \cdot) $ denote a chosen $ G $-invariant inner product on $ \lieg $. Let $ \proots $ denote the
collection of positive roots of $ G $. Let $ \opwc $ denote the open positive Weyl chamber.

Let $ s $ denote the real dimension of $ G $. Let $ l $ denote the real dimension of $ T $.
Let $ N \geq 3 $ be a positive integer.

Let $ \vecxi = (\xi_{1}, \dots, \xi_{N}) $ be an $ N $-tuple of elements in $ \lieg $.
Let $ M = \orb_{\xi_{1}} \times \dots \times \orb_{\xi_{N}} $ be the product of the corresponding
adjoint orbits. Then $ G $ acts on $ M $ through the diagonal adjoint action. This is a Hamiltonian
action, so we have the moment maps $ \mu_{G} $ and $ \mu_{T} $ as before. Let $ M_{0} = \mu_{G}^{-1}(0) $.
Let $ M_{\red} = M_{0} / G $ be the reduced space, which we call an $ N $-fold reduced product of $ G $.

The input $ \vecxi = (\xi_{1}, \dots, \xi_{N}) $ satisfies the following assumptions
 (as in \S 2):
\begin{description}
\item[(A1)] $ \mu_{G}^{-1}(0) \neq \emptyset $ and $ 0 $ is a regular value for $ \mu_{G} $.
\item[(A2)] All $ \xi_{i} $'s lie in $ \opwc $.
\item[(A3)] $ M_{\red}(\vecxi) $ is a smooth manifold.
\end{description}

Given $ \vecw = (w_{1}, \dots, w_{N}) \in W^{N} $, let
\begin{equation}
\vecw \cdot \vecxi = (w_{1} \cdot \xi_{1}, \dots, w_{N} \cdot \xi_{N})
\end{equation}
and
\begin{equation}
\vecw \odot \vecxi = \sum_{i = 1}^{N} w_{i} \cdot \xi_{i},
\end{equation}
just as before. Then, we have
\begin{equation}
M^{T} = \left\lbrace \vecw \cdot \vecxi \  : \  \vecw \in W^{N} \right\rbrace ,
\end{equation}
where $ M^{T} $ denotes the fixed point set of the action of $ T $ on $ M $.
Notice that $ M^{T} $ is discrete and $ \left| M^{T} \right| = \left| W \right| ^{N} $.

Let
\begin{equation}
\varpi(\psi) = \prod_{\gamma \in \proots} \gamma(\psi)
\end{equation}
for all $ \psi \in \liet $.

We have the following general result. 

\begin{theorem}
Let $ G $ be a general semisimple compact connected Lie group and
$ M = \orb_{\xi_{1}} \times \dots \times \orb_{\xi_{N}} $ be the Cartesian product
of $ N \geq 3 $ adjoint orbits of $ G $, where the $ \xi_{i} $'s satisfy the
assumptions (A1), (A2) and (A3). Then,
the symplectic volume of $ M_{\red}(\vecxi) $ is
\begin{equation}
\vol^{\SSS}(M_{\red}(\vecxi)) = \frac{1}{\I^{d/2}} n_{0} C_{G} \res\left( \sum_{\vecw \in W^{N}}
\sgn(\vecw) \frac{e^{\I (\vecw \odot \vecxi, \psi)} [d\psi]}{\varpi^{N-2}(\psi)} \right) ,
\end{equation}
where
\begin{equation}
d = N(s - l) - 2s = (N - 2)s - Nl,
\end{equation}
and $ n_{0} $ is the cardinality of the stabilizer $ \Stab_{G}(p) $ of a generic point $ p $ in $ M_{0} $
and the constant $ C_{G} $ is defined by
\begin{equation}
C_{G} := \frac{(-1)^{n_{+}}}{(2 \pi)^{s - l} \left| W \right| \vol^{\RRR}(T)}.
\end{equation}
Here, $ n_{+} $ is the number of positive roots, in other words, $ n_{+} = (s - l) / 2 $.
\end{theorem}

\begin{proof}
At each fixed point $ F = \vecw \cdot \vecxi \in M^{T} $, the $ T $-equivariant Euler class
of the normal bundle over $ F $ is
\begin{equation}
\ee_{F}(\psi) = \varpi^{N}(\psi).
\end{equation}
The computation  of the residue is similar to that in the case $N=3$ and $G=SU(3)$.
\end{proof}



\bibliographystyle{plain}



\end{document}